\documentclass{article}
\usepackage{amsmath,amsthm,amssymb,amsfonts}
\usepackage{hyperref,lineno}
\usepackage{xcolor}
\usepackage{geometry}

\DeclareSymbolFont{newfont}{OML}{cmm}{m}{it}

\hypersetup{
	colorlinks,
	linkcolor={red!50!black},
	citecolor={blue!50!black},
	urlcolor={blue!80!black}
}

\newtheorem{proposition}{Proposition}
\newtheorem{lemma}{Lemma}[section]
\newtheorem{remark}{Remark}
\newtheorem{corollary}{Corollary}
\newtheorem{theorem}{Theorem}
\newtheorem{example}{Example}
\newtheorem{definition}{Definition}

\title{Maximal Inequalities for Empirical Processes  under General Mixing Conditions} 

\author{Demian Pouzo \\ UC Berkeley \thanks{Dept. of Economics, UC Berkeley. email: dpouzo@berkeley.edu. I would like to thank Michael Jansson for comments, and an anonymous referee and AE for thoughtful comments and pointing out a gap in the initial version. All errors are mine.}}

\date{\today}

\begin{document}

	\maketitle
	\thispagestyle{empty}
	
	\begin{abstract}
    This paper provides a bound for the supremum of sample averages over a class of functions for a general class of mixing stochastic processes with arbitrary mixing rates. Regardless of the speed of mixing, the bound is comprised of a concentration rate and a novel measure of complexity. The speed of mixing, however, affects the former quantity implying a phase transition. Fast mixing leads to the standard root-n concentration rate, while slow mixing leads to a slower concentration rate whose speed depends on the mixing structure. Our findings are applied to obtain new Glivenko-Cantelli type results. 
	\end{abstract}

	\setcounter{page}{1}

\section{Introduction}

Maximal inequalities serve as a fundamental cornerstone in empirical processes theory, playing a pivotal role in deriving crucial results, including but not limited to the functional central limit theorem and strong approximations. 

This paper provides a maximal inequality for $f \mapsto G_{n}[f] : = n^{-1/2} \sum_{i=1}^{n} f(X_{i}) - E_{P}[f(X_{i})]$ of the form
\begin{align}\label{eqn:Max.Ineq}
	  \left \Vert  \sup_{f , f_{0}  \in \mathcal{F}}  |G_{n}[f-f_{0}] |    \right \Vert_{L^{1}(P)} \leq \Gamma
\end{align}
for some quantity $\Gamma$ to be specified below, and where the data $(X_{i})_{i=-\infty}^{\infty}$ is a stationary process drawn from a probability $P$ which satisfies some general mixing conditions to be described below.

There is a large body of work deriving maximal inequalities and their derivatives such as the functional CLT for dependent data under mixing conditions;  cf. \cite{doukhan1987principes,arcones1994central,andrews1994introduction,yu1994rates,DMR1995,billingsley2013convergence}, and \cite{dedecker2002maximal} and \cite{rio2017asymptotic} for reviews. Closest to this paper is the important work in \cite{DMR1995} wherein the authors establish a functional invariance principle in the sense of Donsker for absolutely regular empirical processes, where the constant $\Gamma$ is proportional to a measure of complexity of the class $\mathcal{F}$. To the best of our knowledge, this and all other existing results in this literature are derived under restrictions on the decay of mixing coefficients, e.g. $\sum_{k=1}^{\infty} \beta(k) < \infty$ (\cite{DMR1995}), or $k^{\frac{p}{p-2}} (\log k )^{2 \frac{p-1}{p-2}} \beta(k) = o(1)$  for some $p  >2$ (\cite{arcones1994central}) where $\beta$ is the $\beta$-mixing coefficient used in \cite{DMR1995}.
	
	These results leave open the question of what type of maximal inequality one can obtain in contexts where the mixing coefficients do not satisfy these conditions. Many processes do not satisfy them, either because the data exhibits long-range dependency or long-memory and this feature is modeled using slowly decaying dependence structure, e.g. see \cite{asadi2014stationary} in the context of Markov processes; or the data is described by a so-called infinite memory chain (see \cite{doukhan2008weakly}); or simply because the $\beta$-mixing coefficients decay at a slow polynomial rate (see \cite{CCLJOE10}), or not at all (see \cite{andrews1984non}). More generally, there is the open question of how do the mixing properties affect the concentration rate and the constant $\Gamma$ of the maximal inequality \ref{eqn:Max.Ineq}. This paper provides insights into these questions. 
	
Unfortunately, the approach utilized in \cite{DMR1995} and related papers cannot be applied to establish maximal inequalities when the aforementioned restrictions on the mixing coefficients do not hold. One of the cornerstones of this approach relies on insights that can be traced back to Dudley in the 1960s (\cite{dudley_sizes_1967}) for Gaussian processes. Dudley's work states that the ``natural" topology to measure the complexity of the class of functions $\mathcal{F}$ is related to the variation of the stochastic process. In \cite{DMR1995}, the authors use this insight to construct a ``natural" norm, which turns out to depend on the $\beta$-mixing coefficients. Unfortunately, without restrictions on the mixing coefficients, this approach is not feasible because this norm may not even be well-defined.

In view of this, the current paper proposes a new proof technique which employs a \emph{family} of norms, rather than just one, to measure the complexity of the class $\mathcal{F}$. To accommodate this new feature, we introduce a new measure of complexity inspired by Talagrand's measure (\cite{talagrand2005,talagrand2014}) that allows for a family of norms. 

The family of norms proposed in the proof is linked to the dependence structure, which is captured by suitably chosen mixing coefficients. Papers such as \cite{arcones1994central,DMR1995} use $\beta$-mixing while some other papers use stronger concepts such as $\phi$-mixing (see \cite{ibragimov2012gaussian}). In this paper, however, we use the weaker notions of $\tau$-mixing introduced by \cite{dedecker2004coupling,dedecker2005new} and $\alpha$-mixing (cf. \cite{rio2017asymptotic}). The $\tau$-mixing coefficients not only are  typically weaker than the $\beta$-mixing ones, thereby encompassing a wider class of stochastic processes (see  \cite{dedecker2004coupling,dedecker2005new} and examples below), but they are also adaptive to the size of the class of functions $\mathcal{F}$ --- a property not enjoyed by standard mixing coefficients.

The main result in the paper shows that the $L^{1}$ norm of $\sup_{f, f_{0} \in \mathcal{F}} |G_{n}[f-f_{0}]|$ is bounded (up to constants) by the aforementioned complexity measure. Under IID (or $m$-dependent data), the bound corresponds, up to constants, to that in Talagrand's Generic Chaining approach (\cite{talagrand2014} e.g. Theorem 4.5.16). For general mixing data, when a summability condition, similar to the one in \cite{DMR1995} holds, our bound replicates  (up to constants), and in some cases improves upon, the results in the literature. When the summability restriction does not hold --- i.e., the mixing coefficients do not decay quickly to zero --- a bound of the form \ref{eqn:Max.Ineq} remains valid. However, in this case, the quantity $\Gamma$ is comprised not only of the complexity measure, as in the standard case, but also of a scaling factor that depends on the mixing properties and the sample size. This last novel result implies that the concentration rate is not root-n, as in the standard case, but slower and is a function of the mixing rate.

The remainder of the paper is organized as follows. Section \ref{sec:main.maximal} presents the maximal $L^{1}$-inequality results and Section \ref{sec:proofs} presents the proofs. Some technical lemmas and proofs are relegated to the Appendix.

\section{A Maximal $L^{1}(P)$ Inequality}
\label{sec:main.maximal}

This section aims to establish an upper bound for $\left \Vert \sup_{f , f_{0} \in \mathcal{F}} | G_{n}[f-f_{0}] | \right \Vert_{L^{1}(P)}$, where $\mathcal{F}$ is a class of functions bounded in $L^{2}(P) \cap L^{\infty}$.\footnote{As per \cite{talagrand2014}, the quantity $\left \Vert \sup_{f , f_{0} \in \mathcal{F}} | G_{n}[f-f_{0}] | \right \Vert_{L^{1}(P)}$ is defined as $\sup \{ \left \Vert \sup_{f , f_{0} \in \mathcal{M}}  | G_{n}[f-f_{0}] | \right \Vert_{L^{1}(P)} \colon \mathcal{M} \subseteq \mathcal{F}~\text{finite} \}$ to sidestep measurability issues.} We now outline an informal roadmap for the proof to identify key components and motivate subsequent formal definitions. 

The first step of the proof involves constructing a ``chain" between $f$ and $f_{0}$ denoted as $f - f_{0} = \sum_{k=1}^{\infty} \Delta_{k} f$, where each link, $\Delta_{k}f$, belongs to a class with finite cardinality. As a consequence of this chain, the empirical process is decomposed as $G_{n}[f] = \sum_{k=1}^{\infty} G_{n}[\Delta_{k} f ]$.

For each link of the chain (indexed by $k \in \mathbb{N}$), the second step of the proof couples the empirical process $f \mapsto G_{n}[\Delta_{k} f ]$ with $f \mapsto G^{\ast}_{n}[\Delta_{k} f] : = n^{-1/2} \sum_{i=1}^{n}{ \Delta_{k} f(X^{\ast}_{i}) - E_{P}[ \Delta_{k} f(X^{\ast}) ]}$, where $(X_{i})_{i=-\infty}^{\infty}$ and  $(X^{\ast}_{i})_{i=-\infty}^{\infty}$ have joint probability denoted by $\mathbb{P}$; and the latter process is such that $(U^{\ast}_{2j}(q_{k}))_{j=0}^{\infty}$ form an independent sequence and $(U^{\ast}_{2j+1}(q_{k}))_{j=0}^{\infty}$ form another independent sequence, where each block is given by $U^{\ast}_{i}(q_{k})  : =(X^{\ast}_{q_{k}i+1},...,X^{\ast}_{q_{k}i+q_{k}})$ and has the same distribution as its counterpart in the $(X_{i})_{i=-\infty}^{\infty}$ process --- $q_{k} \in \mathbb{N}$ will be chosen below.\footnote{In this informal presentation it is implicitly assumed that $n/q_{k}$ is an integer; we relax this assumption in the general theory.} Existence of such process is established by the results in \cite{dedecker2006inequalities} which are presented in Appendix \ref{app:coupling} for completeness. Henceforth, we refer to $f \mapsto G^{\ast}_{n}[f]$ as a block-independent empirical process. This coupling yields 
\begin{align*}
 \sup_{f , f_{0} \in \mathcal{F}} | G_{n}[f-f_{0}] |  \leq  \sup_{f \in \mathcal{F}} \sum_{k=1}^{\infty}   | G^{\ast}_{n}[\Delta_{k} f] | + \sup_{f \in \mathcal{F}}  \sum_{k=1}^{\infty}  | G^{\ast}_{n}[\Delta_{k} f] - G_{n}[\Delta_{k} f] | 
\end{align*}
where the $L^{1}$-norm of the second term on the right-hand side is controlled by our measure of dependence defined below. 

The final step of the proof involves bounding the term $\left \Vert \sup_{f \in \mathcal{F}} | \sum_{k=1}^{\infty} G^{\ast}_{n}[ \Delta_{k} f ] | \right \Vert_{L^{1}(P)}$. To achieve this, we employ Talagrand's generic chaining insights (cf. \cite{talagrand2014}). However, it is crucial to adjust the measure of complexity to accommodate the link-specific length of the block, given by $q_{k}$. 
As mentioned in the introduction, the reason why the block length influences the topology employed to gauge complexity comes from the realization that the ``natural" distance used in complexity computations stems from the stochastic process's variability. Due to the dependence structure of $f \mapsto G^{\ast}_{n}[f]$, the variability is measured by $\sigma_{2}(f,q_{k}) : = \sqrt{E_{P} \left[   \left(  q^{-1/2}_{k} \sum_{i=1}^{q_{k}} f(X_{i})  - E_{P} [ f(X) ]  \right)^{2}  \right] }$, where $q_{k}$ is the parameter regulating the length of the blocks. Therefore, the ``natural" notion of distance may differ for each link in the chain. To accommodate for this feature we introduce a novel measure of complexity.

From this brief description of the proof one can see the two key ingredients: First, a measure of dependence, which quantifies the error of approximating the original process with a  block-independent one. Second, a measure of complexity that can accommodate a family of norms. Below we formally define these two quantities, we explain the need to work with a family of norms as opposed to just one as in the rest of the literature, and we then present the main result.

%

\paragraph{Measure of Dependence.} Let $P$ be the probability distribution of the stationary process $(X_{i})_{i=-\infty}^{\infty}$, $P_{X}$ denotes the marginal distribution of $X_{i}$ where $X_{i} \in \mathbb{X}$ for some Polish space $\mathbb{X}$, and $P(.\mid \mathcal{M}_{j}^{l})$ denotes the conditional probability distribution given $\mathcal{M}_{j}^{l}$ for any $-\infty \leq j \leq l \leq \infty$, where $\mathcal{M}_{j}^{l}$ is the $\sigma$-algebra generated by $(X_{j},\ldots,X_{l})$.

The measure that quantifies the dependence structure of the data for any $\mathcal{B} \subseteq cone \mathcal{F}$  is given by\footnote{For any set $A \subseteq L^{2}(P)$, $cone A : = \{ \lambda a \colon a \in A~and~\lambda \geq 0  \}$.}  (cf. \cite{dedecker2006inequalities,dedecker2004coupling,dedecker2005new})
\begin{align}
	& \tau_{\mathcal B}(q)   :  =\max_{1\leq l \leq q} \frac{1}{l} \sup \left\{  	\tau_{d_{l,\mathcal B}}( \mathcal{M}_{-\infty}^{-q}  ; X_{t_{1}},...,X_{t_{l}}) \colon 1\leq t_{1} \leq \ldots \leq t_{l}   \right\},~\forall q \in \mathbb{N},
\end{align}
where the supremum is taken over all future blocks of at most $q$ elements with starting time $t_1 \ge 1$, and 
\begin{align}\label{eqn:tau.defn}
\tau_{d_{l,\mathcal B}}( \mathcal{M}^{-q}_{-\infty}  ; X_{t_{1}},...,X_{t_{l}}) : =  \left \Vert \sup_{g \in \Lambda( \mathbb{X}^{l} , d_{l,\mathcal B}  )}  \left| E_{P}[ g(X_{t_{1}},...,X_{t_{l}})  \mid \mathcal{M}_{-\infty}^{-q}]  -  E_{P}[g(X_{t_{1}},...,X_{t_{l}})] \right| \right \Vert_{L^{1}(P)},
\end{align}
with $\Lambda( \mathbb{X}^{l} , d_{l,\mathcal B}  )$ is the class of Lipschitz (with constant 1) functions over $\mathbb{X}^{l}$ with respect to $(x_{1:l},y_{1:l}) \mapsto d_{l,\mathcal B}(x_{1:l},y_{1:l}) : = \sum_{m=1}^{l} \sup_{f \in  \mathcal{B}} |f(x_{m}) - f(y_{m})|$.\footnote{The element $x_{1:l}$ denotes the vector  $(x_{1},...,x_{l})$ for any $l \in \mathbb{N}$.} 


The quantity $\tau_{\mathcal B}(q)$, introduced in 
\cite{dedecker2006inequalities,dedecker2004coupling,dedecker2005new}, 
measures the largest average per-coordinate dependence between the distant past 
(up to time $-q$) and the future --- measured by any future block of length at most $q$. 
The dependence is quantified in a Kantorovich--Wasserstein sense using 
the class of $1$-Lipschitz functions on $\mathbb X^{l}$ with respect to the 
block metric generated by $\mathcal B$. We defer a more detailed discussion of its properties and its relation to other mixing coefficients, such as $\beta$-mixing, to below. For present purposes, we simply note that this notion of dependence --- unlike $\beta$-mixing --- can be tailored to the function class under consideration, as we do next.

Let
\begin{align}\label{eqn:defn.tautheta}
	q \mapsto \tau(q) : = \tau_{B \mathcal{F} }(q)~and~q \mapsto \theta(q) : = \max\{ 	\tau(q) , \alpha(q)  \},
\end{align}
where $B \mathcal F : = \{  f \in cone (\mathcal{F}-\mathcal{F}) \colon ||f||_{L^{\infty}} \leq 1   \}$, and $\alpha$ is the strong mixing coefficient (e.g. see \cite{rio2017asymptotic})
\begin{align}\label{eqn:alpha.coeff}
	\alpha(q)
	= & \sup_{f \in BL}  \left  \Vert   \int f(x)  P(dx |\mathcal{M}_{-\infty}^{-q} )  -  \int f(x) P_{X}(dx)    \right \Vert_{L^{1}(P)},~\forall q \in \mathbb{N},
\end{align}
where $BL$ is the class of bounded (by 1) measurable functions. Henceforth, we extend $\theta(.)$ to the positive reals as a cadlag function --- to ease the notational burden we still use $\theta$ to denote this extension.

We say the process is $z$-mixing if $\lim_{q \rightarrow \infty} z(q) = 0$ with $z \in \{\theta,\tau,\alpha\}$. The quantity $\theta $ is the relevant quantity for measuring dependence as both $\tau(q)$ and $\alpha(q)$ are used in the proof. The first one controls the error of coupling the original empirical process with a block-independent one; see Lemma \ref{lem:tau-bdd} in Appendix \ref{app:coupling}. On the other hand, as shown in the proof of Lemma \ref{lem:bere} in Appendix \ref{app:main.proof}, the second mixing coefficient, $\alpha$, is used to control the variance of the block-independent empirical process given by  $\sigma^{2}_{2}(f,q) = E_{P} \left[   \left(  q^{-1/2} \sum_{i=1}^{q} f(X_{i})  - E_{P} [ f(X) ]  \right)^{2}  \right] $.

In previous work (cf. \cite{DMR1995,yu1994rates}) the two dependence measures, $\tau$ and $\alpha$, were subsumed by the $\beta$-mixing coefficient, $\beta(q) : = \beta(\mathcal{M}_{-\infty}^{-q}, \mathcal{M}_{0}^{+\infty}) $, where $\mathcal{M}_{0}^{\infty}$ is the $\sigma$-algebra generated by the ``future" $(X_{0},X_{1},...)$ and $\beta$ is the defined as in \cite{VolRoz1959}. The distinction between this $\beta$-mixing coefficient and $\tau$ is twofold. First, $\beta$-mixing compares the distance past $\sigma$-algebra and entire future $\sigma$-algebra, whereas $\tau(q)$ only measures dependence between the distant past $\sigma$-algebra and finite future blocks. Second, and more importantly, $\beta$-mixing is defined through total variation (i.e., supremum over the class $BL$), while $\tau$ relies on a weaker Kantorovich-type metric restricting attention to the class of $1$-Lipschitz functions generated by $B\mathcal F$. So, for any $q \in \mathbb{N}$, $\tau(q)$ is smaller (up to a constant) than $\beta(q)$ --- Lemma \ref{lem:tau.bound.beta} in the Appendix \ref{app:tau.bound} formally shows this.

This class of 1-Lipschitz functions generated by $B\mathcal F$ will be equal to $BL$ when $\mathcal{F}$ is sufficiently rich so that $d_{\mathcal{F}}(x,y)= 1 \{x \ne y \}$;  for instance, if $\mathcal{F}$ is the class of indicators over the half-line.  In many applications, however, it is common for $\mathcal{F}$ to have smoothness restrictions of some sort. These restrictions imply that  the class of $1$-Lipschitz functions generated by $B\mathcal F$ will be much smaller than $BL$,  so that dependence conditions based on $\beta$-mixing may be unnecessarily restrictive. 
 For instance, if functions in $B \mathcal{F}$ are Lipschitz with respect to some distance $\rho$, our $\tau$-coefficient is bounded by that in \cite{dedecker2006inequalities,dedecker2004coupling} as the next lemma shows.

\begin{lemma}\label{lem:tau-bdd.Lip.text}
	For any $q \in \mathbb{N}$, 
	\begin{align*}
		\tau(q) \leq \max_{1\leq l \leq q} \frac{1}{l} \sup \left\{  	\tau_{\rho^{(l)}}( \mathcal{M}^{-q}_{-\infty}  ; X_{t_{1}},...,X_{t_{l}}) \colon 1\leq t_{1} \leq \ldots \leq t_{l}   \right\} = : \tau_{q,\rho}(q)
	\end{align*}
	where, $\rho^{(l)}(x_{1:l},y_{1:l}) : = \sum_{m=1}^{l} \rho(x_{m},y_{m})$. 
\end{lemma}

\begin{proof}
	See Appendix \ref{app:Lipschitz} with $\mathcal B$ equal to $B \mathcal F$. 
\end{proof}
This connection is useful because the literature on $\tau$-dependence (cf \cite{dedecker2004coupling,dedecker2005new}) provides many examples of processes 
that are $\tau$-mixing even though they fail to be $\beta$-mixing --- even cases where the $\tau$-mixing coefficients vanish at exponential rate; one such case is discussed in  Example \ref{exa:AR.discrete} below.

\paragraph{Measure of Complexity. }  We present a new measure of complexity which we dub Talagrand's complexity measure as it is inspired by Talagrand's Generic Chaining theory (see \cite{talagrand1996, talagrand2005, talagrand2014}).


For any $r > 0$ and $\mathcal{B} \subseteq \mathcal{F}$, we say $\mathcal{T}^{\infty} : = (\mathcal{T}_{l})_{l \in \mathbb{N}_{0}}$ is an admissible partition sequence of order $r$ for $\mathcal{B}$, if $\mathcal{T}^{\infty}$ is an increasing sequence of partitions  of  $\mathcal{B}$ with $card \mathcal{T}_{l} \leq 2^{2^{l/r}}$ and $card \mathcal{T}_{0} = 1$.\footnote{By increasing we mean that any set in $\mathcal{T}_{l+1}$ is included in a set in $\mathcal{T}_{l}$.} Let $\mathbf{T}_{r}(\mathcal{B})$ denote the set of all admissible partition sequences of order $r$ for $\mathcal{B}$. For any $f \in \mathcal{B}$ and $l \in \mathbb{N}_{0}$, let $T(f,\mathcal{T}_{l})$ be the (only) set in $\mathcal{T}_{l}$ containing $f$ and let $x \mapsto D(f,\mathcal{T}_{l})(x) : = \sup_{f_{1},f_{2} \in T(f,\mathcal{T}_{l})} |f_{1}(x) - f_{2}(x)|$ be the ``diameter" of such set. 
	
	Talagrand's complexity measure is defined by
	
\begin{definition}[Talagrand's complexity measure]\label{def:MoC}
	For a set $\mathcal{B} \subseteq \mathcal{F}$, order   $r > 0$, index $p>0$, and a family of quasi-norms, $\mathbf{d} : = (d_{l})_{l \in \mathbb{N}_{0}}$, let 
	\begin{align}
		\gamma_{p,r}(\mathcal{B},\mathbf{d}) =     \inf_{\mathcal{T}^{\infty} \in \mathbf{T}_{r}(\mathcal{B})} \sup_{f \in \mathcal{B}}  \sqrt{2} \sum_{l=0}^{\infty} 2^{l/p}   d_{l}(D(f,\mathcal{T}_{l}))
	\end{align}
	be the Talagrand's complexity Measure of set $\mathcal{B}$ under the family $\textbf{d}$ (indexed by $p,r$). 
\end{definition}

The main difference with the expression in definition 2.2.19 and the expression in p. 32  in \cite{talagrand2014} is the usage of a different (quasi)-norm, $d_{l}$,  for each different partition $\mathcal{T}_{l}$.  For reference, when specializing to one norm (as opposed to a family), this definition with $r=1$ is analogous to that in \cite{talagrand2014}. The order $r$ (and the index $p$) provide flexibility and allow us to majorize the measure of complexity by the standard Dudley's metric entropy --- see Lemma \ref{lem:bound.Dudley} in Appendix \ref{app:bound.Dudley}.

\paragraph{Notion of distance.} We now introduce the notion of distance, captured by a family of norms. The ``natural" notion of distance to compute the complexity arises from the Berenstein inequality (see Lemma \ref{lem:bere} below), which works with the boundedness and the variability of the stochastic process. As argued above the relevant process is the block empirical process $f \mapsto G^{\ast}_{n}[f]$, and due to its dependence structure the variability is given by $\sigma_{2}(f,q)$, where $q$ is the parameter regulating the length of the blocks. Motivated by this observation and Theorem 1.1 in \cite{rio2017asymptotic}, we introduce the following norms
\begin{align}
	||f||^{2}_{2,q} : = \int_{0}^{1} \mu_{q}(u) Q^{2}_{f}(u) du~and~||f||_{\infty,q} : = q ||f||_{L^{\infty}} : = q \sup_{x} |f(x)|,~\forall f \in \mathcal{F},\forall q \in \mathbb{N},
\end{align}
where 
\begin{align}
	\label{eqn:muq}
	u \mapsto \mu_{q}(u) = \sum_{i=0}^{q} 1_{\{ u \leq \alpha (i)  \}} \in \{0,...,1+q\},~\forall q \in \mathbb{N} 
\end{align}
with $u \mapsto Q_{f}(u) : =  \inf \{ s \mid H_{f}(s) \leq u  \}$ being the quantile function of $|f|$ with  $s \mapsto H_{f}(s) : = P(|f(X)| > s)$. 

%
%

\paragraph{Choice of block length.}	Most of the literature derives maximal inequalities of the type studied here under the assumption that $\beta$ mixing coefficients decay ``fast enough" to zero ---  $\sum_{i=1}^{\infty} \beta(i) < \infty$ or stronger. Under this assumption, for any $q$, $||f||_{2,q}$ is majorized (up to constants) by $ \sqrt{\int_{0}^{1} \beta^{-1}(u) Q^{2}_{f}(u) du}$, which is a well-defined norm and can be combined with known metric entropies  --- to our knowledge, such approach was first proposed in \cite{DMR1995} for constructing bracketing entropies using $\beta$-mixing coefficient.  If $\sum_{i=1}^{\infty} \beta(i)$ is not finite, however, the above proposal becomes infeasible since $\beta^{-1}$ is not integrable and thus the previous norm is not even be well-defined. Consequently, the strategy of proof proposed by \cite{DMR1995} and related papers cannot be followed. 

In order to sidestep this issue, we use a different strategy of proof which relies on directly using $||.||_{2,q}$, as it is always well-defined for any finite $q$, and working with a sequence of block lengths, $(\mathbf{q}_{n}(k))_{k=0}^{\infty}$, 
given by
\begin{align}\label{eqn:qk}
	\mathbf{q}_{n}(k) = & \min \{  s \in \mathcal{Q}_{n} \mid \tau(s)  n \leq s 2^{k/2}  \},~\forall (n , k) \in \mathbb{N} \times \mathbb{N}_{0},
\end{align}
where $\mathcal{Q}_{n} :  =  \{ m \in \mathbb{N} \colon m \leq n   \}$.

\paragraph{Main Result.} We are now in position to state the maximal inequality result.

\begin{theorem}\label{thm:equi}
	For any orders $a,b \geq 1$, there exists a constant $\mathbb{L}_{a,b}$, such that for any $n \in \mathbb{N}$,\footnote{The constant $\mathbb{L}_{a,b} $ depends on the orders $a,b$ but is otherwise universal; in particular it does not depend on $P,\mathcal{F},$ or $n$. Its exact expression can be found in the proof.} 
	\begin{align}
		\left \Vert   \sup_{f,f_{0} \in \mathcal{F}}  | G_{n}[f - f_{0}] |    \right \Vert_{L^{1}(P)}  \leq \mathbb{L}_{a,b} \left( \frac{\gamma_{1,b}( \mathcal{F} , ||.||_{\infty,\mathbf{q}_{n}}  ) }{\sqrt{n}}  + \gamma_{2,a}( \mathcal{F} , ||.||_{2,\mathbf{q}_{n}}  ) \right).
	\end{align}
\end{theorem}	

\begin{proof}
	See Section \ref{sec:proof.main}. 
\end{proof}

Unlike classical maximal inequalities, which typically assume independence or fast-mixing dependence (e.g., summable $\beta$-mixing coefficients), this result is valid for arbitrarily slow mixing rates. The bound consists of two complexity measures $\gamma_{1,b}$ and $\gamma_{2,a}$, the first one involving the sup norm and the second one involving an $L^{2}$-norm. This feature is analogous to the results in  Theorem 4.5.16 in \cite{talagrand2014} for IID data.  The key difference lies on the norms being defined relative to a ``mixing-adjusted" family of norms, making the bound adaptive to different dependence structures.\footnote{If the quantities in the RHS are infinite the inequality is trivially true, so, henceforth we assume that are finite. We provide examples wherein this holds below and in Remark \ref{rem:bound.gamma1} in Appendix \ref{app:bound.Dudley}.}

In what follows, we present a few remarks and corollaries designed to shed light on the scope of the theorem, how to use it, its relationship with existing literature, and the role of the dependence structure. 

\subsection{Interpretation and Implications of Theorem \ref{thm:equi}} 

Henceforth, let $x \mapsto \bar{\tau}(x) : = \tau(x)/x$, and $\bar{\tau}^{-1}$ its generalized inverse. We adopt the following convention: $\mathbb{L}$ represents an universal constant and $\mathbb{L}_{x}$ represents an universal constant with the exception that can depend on a parameter ``x" --- e.g., in Theorem \ref{thm:equi}, $\mathbb{L}_{a,b}$ is an universal constant that only depends on $(a,b)$. The constants $\mathbb{L}_{x}$ can take different values at different instances, and are used to simplify the exposition. The derivations in the proofs contain the exact values behind these constants.

\subsubsection{Mixing structure and Geometric structure}

We begin with a refinement of the previous theorem that provides an upper bound that separates the dependence structure from the geometric one. Although this bound is looser than the one in the theorem, it is perhaps more practical and easier to compute in applications. Henceforth, for any $q \in \mathbb{N}$ and any $r \in [1,+\infty]$, let $ \Theta_{r}(q)  : = \sqrt{ 1 + ( \int_{0}^{1}  | \min\{ \alpha^{-1}(u), q  \} |^{\frac{r}{r-1}} du  )^{\frac{r-1}{r}}  } $ for $r \in (1,+\infty)$, $\Theta_{1}(q)  : = \sqrt{1+q} $ for $r=1$, and $\Theta_{\infty}(q)  : = \sqrt{1 +  \int_{0}^{1}  | \min\{ \alpha^{-1}(u), q  \} | du    } $ for $r=+\infty$.

\begin{corollary}\label{cor:bound.separation}	
	For any orders $a,b \geq 1$, any $r \in [1,+\infty]$ and any $n \in \mathbb{N}$,
	\begin{align}
		\left \Vert   \sup_{f,f_{0} \in \mathcal{F}}  | G_{n}[f - f_{0}] |    \right \Vert_{L^{1}(P)}  \leq \mathbb{L}_{a,b} \left( \frac{\bar{\tau}^{-1}(n^{-1})}{\sqrt{n}}   \gamma_{1,b}( \mathcal{F} , ||.||_{L^{\infty}}  )  +  \Theta_{r} (\bar{\tau}^{-1}(n^{-1}))      \gamma_{2,a}( \mathcal{F} , ||.||_{L^{2r}(P)}  ) \right),
	\end{align}

	
\end{corollary}

\begin{proof}
	See Section \ref{sec:bound.separation}. 
\end{proof}

This result is obtained by showing that, for any $q \in \mathbb{N}$,
\begin{align}\label{eqn:bound.separation}
 \gamma_{1,b}( \mathcal{F} , ||.||_{\infty,q}  )\leq q   \gamma_{1,b}( \mathcal{F} , ||.||_{L^{\infty}}  ) ~and~ \gamma_{2,a}( \mathcal{F} , ||.||_{2,q}  )	  \leq \Theta_{r} (q)      \gamma_{2,a}( \mathcal{F} , ||.||_{L^{2r}(P)}  ),
\end{align}
and it provides insight into the role of the mixing structure and the complexity of $\mathcal{F}$. The influence of the mixing structure is reflected in the scaling factors $\frac{\bar{\tau}^{-1}(1/n)}{\sqrt{n}}$ and $\Theta_{r}(\bar{\tau}^{-1}(1/n))$, which depend on the behavior of $\tau$ through $\bar{\tau}^{-1}$ and on $\alpha$ through $\Theta_{r}$. The complexity of $\mathcal{F}$ is captured by $\gamma_{1,b}$ and $\gamma_{2,a}$. These quantities, introduced in Talagrand's work, are independent of the dependence structure and benefit from numerous established upper bounds that we can leverage.

For instance, by Lemma \ref{lem:bound.Dudley} in Appendix \ref{app:bound.Dudley},
\begin{align}  \label{eqn:gamma.entropy.bound}
	&~\gamma_{a,b}(\mathcal{F},||\cdot||_{L^{r}(P)}) \leq \sum_{l=0}^{\infty} 2^{l/a} e_{l,b}(\mathcal{F},||\cdot||_{L^{r}(P)}),\\ \notag
	&~and~if~b\leq 0.5a,~it~follows~that\\ \notag
	&\gamma_{a,b}(\mathcal{F},||\cdot||_{L^{r}(P)}) \leq \mathbb{L} \int_{0}^{Diam(\mathcal{F},||\cdot||_{L^{r}(P)})} \sqrt{\log N(u,\mathcal{F},||\cdot||_{L^{r}(P)}) }du
\end{align}
for any $r \in [1,\infty]$, where $e_{l,b}(\mathcal{F},||\cdot||_{L^{r}(P)}) := \inf_{S \subseteq \mathcal{F}, card S \leq 2^{2^{l/b}}} \sup_{t \in \mathcal{F}} \min_{s \in S} ||t - s||_{L^{r}(P)}$ denotes the entropy number, and $N(u,\mathcal{F},||\cdot||_{L^{r}(P)})$ is the covering number for radius $u$.\footnote{For convenience, the entropy number is defined using a cardinality of $2^{2^{l/b}}$ instead of the conventional $2^{l/b}$.} Expression \ref{eqn:gamma.entropy.bound} is particularly useful as it relies on a well-known and extensively studied quantity: the entropy numbers. It can also be improved; for instance, when $\mathcal{F}$ exhibits certain convexity properties, this bound can be further improved by utilizing the results from \cite{VanHandel2018,VanHandel2018b}.

\subsubsection{The $m$-dependent case}

This case is presented just to illustrate our bound relative to the literature and the effect of dependence in maximal inequalities.  Suppose $(X_{i})_{i}$ is such that, for any $i \in \mathbb{N}$, $X_{i} = \Phi(Z_{i+1},...,Z_{i+m})$ where $(Z_{i})_{i}$ is an IID process. This process is $m$-dependent and so $\theta(q) = 0$ for any $q \geq m$. Therefore, by choosing $\mathbf{q}=m$, it follows that $f \mapsto ||f||^{2}_{2,m} = \sum_{i=0}^{m} \int_{0}^{0.5\theta(i)} Q^{2}_{f}(u) du \leq  m \int_{0}^{1} Q^{2}_{f}(u) du = m ||f||^{2}_{L^{2}(P)}$. Hence, by Theorem \ref{thm:equi}
\begin{align}\label{eqn:bound.mdepen}
	\left \Vert  \sup_{f , f_{0}  \in \mathcal{F}} | G_{n}[f-f_{0}] |  \right \Vert_{L^{1}(P)} \leq  \mathbb{L}_{a,b} \sqrt{m} \left( \sqrt{\frac{m}{n}} \gamma_{1,b}(\mathcal{F} , ||.||_{L^{\infty}}) + \gamma_{2,a}(\mathcal{F} , ||.||_{L^{2}(P)}) \right).
\end{align}

For the IID case of $m=1$, this bound coincides (up to constants) with known bounds in the literature; cf. Proposition 9.2 in \cite{talagrand2014}. 

For a general $m>1$, however, the IID bound is scaled up by a factor of $m$ thereby reflecting the fact that in this type of processes, without additional knowledge on the function $\Phi$, is as if the sample size is $n/m$ as opposed to $n$. 
This matters for applications as it is common practice in statistics to employ an asymptotic approach, wherein $n$ divergences and constants are ignored. This practice will lead to the conclusion that the bound coincides with the IID one, as $m$ is merely constant, but the theorem illustrates that, even with this simple dependence structure, such a conclusion can be misleading. For instance, even for moderate dependence --- e.g., $m$ equal 4 or 5 --- the \emph{actual} bound corresponds to more than twice the IID bound, or alternatively, corresponds to an IID bound but with less than 25 percent of the sample size.

This feature, whereby the dependence structure results in a scaling of the IID bound, is canonical, appearing in more general processes as shown below. Moreover, the nature of the scaling --- whether is constant or diverging with the sample size --- depends on how fast mixing occurs.

\subsubsection{The $\sum_{q=1}^{\infty} \theta(q) < \infty$ case}

In this case, the norm $f \mapsto |||f|||_{2,\theta } : = \sqrt{\int_{0}^{1} (1  + \theta^{-1}(u) ) Q^{2}_{f}(u) du }$ is well-defined as $\theta$ is integrable (cf. \cite{DMR1995}). Moreover, it turns out that $||f||_{2,q} \leq |||f|||_{2,\theta } $ (this is formally shown in the proof of the proposition below). Hence, Theorem \ref{thm:equi} implies the following result. 

\begin{proposition}\label{pro:bound.tau.l1}
	Suppose $\sum_{q=1}^{\infty} \theta(q) < \infty$ and 	$\gamma_{1,b}( \mathcal{F} , ||.||_{L^{\infty}} )  < \infty$ for some order $b \geq 1$. Then, for sufficiently large $n$,
	\begin{align*}
		\left \Vert  \sup_{f , f_{0}  \in \mathcal{F}} | G_{n}[f-f_{0}] |  \right \Vert_{L^{1}(P)} \leq  \mathbb{L}_{b}  \gamma_{2,1}( \mathcal{F} , |||.|||_{2,\theta}  ).
	\end{align*}
\end{proposition}

\begin{proof}[Proof of Proposition \ref{pro:bound.tau.l1}]
	The claim follows by majorizing the RHS of the inequality in Theorem \ref{thm:equi}, i.e., 	for sufficiently large $n$,
	\begin{align}\label{eqn:bound.tau.l1-1}
	\frac{\gamma_{1,b}( \mathcal{F} , ||.||_{\infty,\mathbf{q}_{n}}  )}{\sqrt{n}}  +	\gamma_{2,a}( \mathcal{F} , ||.||_{2,\mathbf{q}_{n}}  ) \leq \mathbb{L}  \gamma_{2,a}( \mathcal{F} , |||.|||_{2,\theta}  ).
	\end{align}

	We first show that for any order $a \geq 1$, $\gamma_{2,a}( \mathcal{F} , ||.||_{2,\mathbf{q}_{n}}  ) $ can majorized (up to constants) by $\gamma_{2,a}( \mathcal{F} , |||.|||_{2,\theta}  )$.  This follows because, for any $q \in \mathbb{N}$, $||f||^{2}_{2,q} \leq |||f|||^{2}_{2,\theta}   = \int_{0}^{1} (1+\theta^{-1}(u)) Q^{2}_{f}(u) du $ for any $f \in \mathcal{F}$. This last inequality holds because $u \mapsto \mu_{q}(u) \leq  \theta^{-1}(u) +1 $ for any $q \in \mathbb{N}$ (see Lemma \ref{lem:mu.bound} in Appendix \ref{app:supp.lemmas}) and, under $\sum_{q=1}^{\infty} \theta(q) < \infty$, $\theta^{-1}$ is integrable. Thus, for any admissible partition sequence of order $a$ for $\mathcal{F}$,
	\begin{align*}
	\sup_{f \in  \mathcal{F}} 	\sum_{l=0}^{\infty} 2^{l/2} ||D(f,\mathcal{T}_{l})||_{2,\mathbf{q}_{n}(l)}\leq \mathbb{L} 	\sup_{f \in  \mathcal{F}} 	\sum_{l=0}^{\infty} 2^{l/2} |||D(f,\mathcal{T}_{l})|||_{2,\theta},
	\end{align*}
so the desired result follows from optimizing over the partition sequence.

	Second, we show that $\limsup_{n\to\infty}	\frac{\gamma_{1,b}( \mathcal{F} , ||.||_{\infty,\mathbf{q}_{n}}  )}{\sqrt{n}}  = 0 $. Since $\mathbf{q}_{n}(0) \geq \mathbf{q}_{n}(l)$ for all $l \in \mathbb{N}_{0}$, it readily follows that $\gamma_{1,b}( \mathcal{F} , ||.||_{\infty,\mathbf{q}_{n}}  ) \leq \mathbf{q}_{n}(0) \gamma_{1,b}( \mathcal{F} , ||.||_{\infty}  )$. Thus, it suffices to show $\limsup_{n\to\infty} \mathbf{q}_{n}(0)/\sqrt{n}  =0$. But this follows from Lemma \ref{lem:tau.bounded.qn.rate} in Appendix \ref{app:supp.lemmas}. 
	
	Hence, equation \ref{eqn:bound.tau.l1-1} holds and by Theorem \ref{thm:equi} with $a=1$ and $\mathbb{L}_b : = \mathbb{L}_{1,b}$ the desired results holds.
\end{proof}

This case of ``fast mixing" has been studied in the literature (cf. \cite{DMR1995} Theorem 1), but Proposition \ref{pro:bound.tau.l1} offer extensions in two directions.  First, by virtue of the generic chaining theory, it provides a tighter bound than existing results which use Dudley's metric entropy or Ossiander's bracketing entropy. Second, and perhaps more importantly, it imposes weaker restrictions on the mixing structure of the data because it relies on the $\theta$-mixing coefficients as oppose to the $\beta$-mixing ones, thus extending  maximal inequalities to a wider class of stochastic processes --- Example \ref{exa:AR.discrete} illustrates this point.

\begin{example}\label{exa:AR.discrete}
	Suppose the process $(X_{i})_{i}$ follows the following auto-regressive model $X_{i} = H(X_{i-1}) + \zeta_{i}$ with $\zeta_{i}$ IID with bounded PDF, $H$ being $\kappa$-Lipschitz, and $\zeta_{0}$ having finite moments. Suppose that $\mathcal{F}$ is the class of monotone functions in $\mathcal{C}^{\eta}([0,1])$ --- the H\"{o}lder class with $\eta$ smoothness, where $\eta > 2$.\footnote{See Remark \ref{rem:bound.gamma1} in Appendix \ref{app:bound.Dudley} for a formal description of this space and others.} 
	
	Since $\mathcal{F} \subseteq \mathcal{C}^{\eta}([0,1])$, $\mathcal{F}$ is contained in the class of Lipschitz functions under the Euclidean metric. Thus, $\Lambda(\mathcal{F})$ is also contained in this class, so $\tau(q)$ is majorized (possibly up to constants) by $\tau_{||.||,q}(q)$, the $\tau$-coefficient in \cite{dedecker2006inequalities,dedecker2005new} (see also Lemma \ref{lem:tau-bdd.Lip.text}). Thus $\tau(q) = O(e^{q \log \kappa})$ (cf.   \cite{dedecker2006inequalities,dedecker2005new}).  
	
	Moreover, by expression \ref{eqn:gamma.entropy.bound}, $\gamma_{2,1}( \mathcal{F}, |||.|||_{2,\theta}) \leq \int_{0}^{M} \sqrt{\log N(u,\mathcal{F},|||.|||_{2,\theta})}du$ where $M : = \sup_{f \in \mathcal{F}} ||f||_{L^{\infty}}$. By the fact that $e_{l,1}( \mathcal{C}^{\eta}([0,1]), ||.||_{L^{\infty}}) \leq \mathbb{L} 2^{- l \eta}$ (see Remark \ref{rem:bound.gamma1} in Appendix \ref{app:bound.Dudley}), $\gamma_{1,1}( \mathcal{F}, ||.||_{L^{\infty}}) < \infty$.

	Hence, by Proposition \ref{pro:bound.tau.l1} 
	\begin{align}\label{eqn:exa.Lip-1}
		\left \Vert  \sup_{f , f_{0}  \in \mathcal{F}} | G_{n}[f-f_{0}] |  \right \Vert_{L^{1}(P)} \leq \mathbb{L} \int_{0}^{M} \sqrt{\log N(u,\mathcal{C}^{\eta}([0,1]),|||.|||_{2,\theta})}du. 
	\end{align}

	By H\"{o}lder inequality, $|||.|||_{2,\theta} \leq \sqrt{||\theta^{-1}||_{L^{1}([0,1],Leb)}} ||.||_{L^{\infty}}$, and $||\theta^{-1}||_{L^{1}([0,1],Leb)} <\infty$ when $\kappa < 1$. So, the RHS in the previous display is bounded (up to constants) by   $\int_{0}^{1} \sqrt{\log N(u,\mathcal{C}^{\eta}([0,1]),||.||_{L^{\infty}})}du$ which is finite provided $\eta > 2$ (see \cite{VdV-W1996} Ch 2.7). 
	
An immediate implication of \eqref{eqn:exa.Lip-1} is a Glivenko–Cantelli type result at rate $n^{-1/2}$ for “fast” $\theta$-mixing autoregressive processes. These processes are not necessarily $\beta$-mixing --- e.g., $\zeta_{i}$ has a discrete distribution ---  and thus standard results in the literature cannot be applied to obtain this result. 
   $\triangle$ 
\end{example}

\subsubsection{Maximal Inequality for slow mixing processes}

In this section we discuss the case of ``slow mixing" --- defined by $\min\{ \sum_{q} \alpha(q) , \sum_{q} \tau(q)\}  = \infty$ ---, which, to our knowledge, has not been studied before in the literature. Intuitively, in this case there is a ``phase-transition" because the scaling parameter due to the mixing structure is no longer constant, rather, it diverges with the sample size due to the slow mixing. That is, while the bound preserves the standard structure of a scaling factor and a measure of complexity, now the scaling factor is divergent with the sample size. 

To formalize this intuition, one can use Corollary \ref{cor:bound.separation}, but  a more succinct refinement can be obtained for $r=1$ if $\mathcal{F}$ is separable under $L^{2}(P)$ norm (a rather mild condition for many applications). Henceforth, let $\varphi : = (\varphi_{l})_{l}$ be the orthonormal basis, so that any $f$ in $\mathcal{F}$ can be cast as $\langle \pi(f) ,  \varphi \rangle_{\ell^{2}}$ for some $\pi(f) \in \Pi(\mathcal{F}) \subseteq \ell^{2}$ where $\Pi(\mathcal{F})$ summarizes additional summability restrictions on the coefficients that $\mathcal{F}$ may induce through its smoothness properties.

\begin{proposition}\label{pro:bound.slow.gaussian}
	Suppose $\gamma_{1,b}( \mathcal{F} , ||.||_{L^{\infty}} )  < \infty$ for some order $b\geq 1$. Then, for sufficiently large $n$, 
	\begin{align}
		\left \Vert   \sup_{f,f_{0} \in \mathcal{F}}  | G_{n}[f - f_{0}] |    \right \Vert_{L^{1}(P)}  \leq \mathbb{L}_{b} \sqrt{\bar{\tau}^{-1}(1/n)} E \left[ \sup_{\pi \in  \Pi(\mathcal{F}) } \langle \pi , \zeta \rangle_{\ell^{2}}  \right] 
	\end{align}
	where $\zeta : = (\zeta_{l})_{l}$ are independent Gaussian random variables.
\end{proposition}

\begin{proof}[Proof of Proposition \ref{pro:bound.slow.gaussian}]
By 	Corollary \ref{cor:bound.separation} with $r=1$, 
	\begin{align*}
	\left \Vert   \sup_{f,f_{0} \in \mathcal{F}}  | G_{n}[f - f_{0}] |    \right \Vert_{L^{1}(P)}  \leq \mathbb{L}_{a,b}  \left( \frac{\mathbf{q}_{n}(0)}{\sqrt{n}}  \gamma_{1,b}( \mathcal{F} , ||.||_{L^{\infty}} )  + \sqrt{ 1 + \mathbf{q}_{n}(0)} \gamma_{2,a}( \mathcal{F} , ||.||_{L^{2}(P)} ) \right). 
\end{align*}
	
	We now show that $\mathbf{q}_{n}(0) = o(n)$. This claim follows because, by definition \ref{eqn:qk}, $\mathbf{q}_{n}(0)$ is either $\mathbf{q}_{n}(0)  = 1$ or $\bar{\tau}(\mathbf{q}_{n}(0)-1) \geq 1/n  \iff n \tau(\mathbf{q}_{n}(0)-1) + 1 \geq \mathbf{q}_{n}(0)$  and since $\lim_{q \rightarrow \infty} \tau(q) = 0$, it follows that $\mathbf{q}_{n}(0) = o(n)$. 
	This and finiteness of $ \gamma_{1,b}( \mathcal{F} , ||.||_{L^{\infty}(P)} )  $ imply that for sufficiently large $n$, the RHS of the previous display is dominated by $\mathbb{L}_{a,b}  \sqrt{1 + \mathbf{q}_{n}(0)} \gamma_{2,a}( \mathcal{F} , ||.||_{L^{2}(P)} ) $. By Lemma \ref{lem:q0.bound}, $1+\mathbf{q}_{n}(0) \leq 2 \bar{\tau}^{-1}(1/n)$ for sufficiently large $n$, thus
	\begin{align*}
	\left \Vert   \sup_{f,f_{0} \in \mathcal{F}}  | G_{n}[f - f_{0}] |    \right \Vert_{L^{1}(P)}  \leq \mathbb{L}_{a,b}  \sqrt{ \bar{\tau}^{-1}(1/n)} \gamma_{2,a}( \mathcal{F} , ||.||_{L^{2}(P)} ).
\end{align*}

	Hence, to obtain the desired result it suffices to bound $\gamma_{2,1}( \mathcal{F} , ||.||_{L^{2}(P)} )$ by $E[\sup_{\pi \in  \Pi(\mathcal{F}) } \langle \pi , \zeta \rangle_{\ell^{2}}  ] $. To show this, we invoke the majorizing measure theorem (MMT; \cite{talagrand2014} Theorem 2.4.1) that yields $\gamma_{2}(T ,d ) \leq \mathbb{L}E[\sup_{t \in T} X_{t}] $ where $t \mapsto X_{t}$ is a Gaussian process, and $d(a,b)^{2} : = E[|X_{a} - X_{b}|^{2}]$.
	
	We now specialize this result to $T = \mathcal{F}$ and $(a,b) \mapsto d(a,b) = ||a-b||_{L^{2}(P)}$. $\mathcal{F}$ is separable with  orthonormal basis given by $(\varphi_{l})_{l}$, any $f \in \mathcal{F}$ is fully characterize by its ``fourier" coefficients $\pi(f) \in \Pi(\mathcal{F})$. We claim that $X_{f} = \sum_{l=1}^{\infty} \zeta_{l} \pi_{l}(f)$ for any $f \in \mathcal{F}$. To show this, we invoke the Karhunen–Lo\'{e}ve (KL) representation theorem. Since $E[|X_{f} - X_{g}|^{2}] = ||f-g||^{2}_{L^{2}(P)} = \sum_{l} ( \pi_{l}(f) - \pi_{l}(g) )^{2}$, the covariance kernel is such that $C(f,g) = \langle f ,g \rangle_{L^{2}(P)} = \sum_{l} \pi_{l}(f) \pi_{l}(g)$, and according to this representation of the covariance kernel, 
	\begin{align*}
		X_{f} = \sum_{l=1}^{\infty} \zeta_{l} \pi_{l}(f),~\forall f \in \mathcal{F},
	\end{align*}
	and so $E[\sup_{t \in T} X_{t}]  = E[\sup_{\pi \in \Pi(\mathcal{F})  } \langle \pi , \zeta \rangle_{\ell^{2}}  ] $. 
\end{proof}

At first glance, Proposition \ref{pro:bound.slow.gaussian} might seem a surprising result as it majorizes the expectation of the supremum of an slowly $\theta$-mixing empirical process with a seemingly unrelated quantity: the expectation of the supremum of a Gaussian process, scaled by $\bar{\tau}^{-1}(1/n)$. However, the result stems from (a) a careful control of the mixing structure and geometry of $\mathcal{F}$, presented in Corollary \ref{cor:bound.separation}, and (b) the celebrated majorizing measure theorem  (e.g. \cite{talagrand2014} Theorem 2.4.1).\footnote{The results and technique in Proposition  \ref{pro:bound.slow.gaussian}  is not confined to slow $\theta$-mixing processes, it can also be obtained for the fast $\theta$-mixing case.}

This result presents a relatively easy to use tool to bound $	\left \Vert   \sup_{f,f_{0} \in \mathcal{F}}  | G_{n}[f - f_{0}] |    \right \Vert_{L^{1}(P)}  $: Computing $\bar{\tau}^{-1}(1/n)$ for a given $\tau$ is trivial, and bounds for the supremum of the Gaussian process has been extensively studied in the literature; e.g. One can use Dudley's metric entropy, or the results in \cite{talagrand2014} or \cite{VanHandel2018}, or simply bounded using H\"{o}lder inequality. In addition, it has the interesting property that it does not depend on the $\alpha$-mixing coefficients, only on the $\tau$-mixing ones, thus, presenting itself as a useful bound even if the process is not $\alpha$-mixing.  
 This simplicity, however, may come at the cost of yielding slower rates of convergence than those implied by Corollary \ref{cor:bound.separation}. 

The next example illustrates these points and also establish new Glivenko-Cantelli type results over Sobolev classes under slow mixing processes.  The key takeaway is that, in this case, the concentration rate is slower than $\sqrt{n}$, revealing a fundamental limitation of classical methods in this setting.

\begin{example}
	Suppose $\tau(q) = \alpha(q)= q^{-m}$ for some $m \in (0,1]$ --- an example of a process with this mixing structure is presented in application 2 in \cite{dedecker2004coupling} and given by an autoregressive process $X_{t+1} = f(X_{t}) + \zeta_{t}$ where $m$ depends on features of $f$ and integrability of $\zeta$.

	 Suppose $\mathcal{F}$ is contained in a Sobolev space with periodic domain, $\mathbb{W}^{s}_{p}(\mathbb{T})$ for some $s>1+1/p$ and $p \geq 1$. It is well-known that $\mathcal{F}$ satisfies the condition of the proposition with $\pi(f)$ being the Fourier coefficients of $f$ and $\Pi(\mathcal{F})  = \{ \pi \in \ell^{p} \colon \sum_{l} (1+|l|^{2})^{sp/2} \pi^{p}_{l} < \infty  \}$. By expression \ref{eqn:gamma.entropy.bound} and  Remark \ref{rem:bound.gamma1} in Appendix \ref{app:bound.Dudley}, $\gamma_{1,1}(\mathbb{W}^{s}_{p}(\mathbb{T}), ||.||_{L^{\infty}}) \leq \sum_{l} 2^{l-ls} < \infty$. 
	 
	 Given our choice of $\tau$, it readily follows that $\bar{\tau}^{-1}(1/n)  = O(  n^{1/(1+m)}   ) $.\footnote{The symbol $O$ stands for bounded (at a particular rate) and the analogous symbol $O_{P}$ stands for bounded in probability $P$ (for a particular rate).}  	Hence,  Proposition \ref{pro:bound.slow.gaussian} implies
		\begin{align*}
		\left \Vert   \sup_{f,f_{0} \in \mathcal{F}}  | G_{n}[f - f_{0}] |    \right \Vert_{L^{1}(P)}  \leq \mathbb{L}_{1} n^{\frac{1}{2(m+1)}}  E \left[ \sup_{\pi \in  \Pi(\mathcal{F}) } \langle \pi , \zeta \rangle_{\ell^{2}}  \right] .
	\end{align*}
By H\"{o}lder inequality and the condition on $\mathcal{F}$,
		\begin{align*}
	\left \Vert   \sup_{f,f_{0} \in \mathcal{F}}  | G_{n}[f - f_{0}] |    \right \Vert_{L^{1}(P)}  \leq \mathbb{L}_{1} n^{\frac{1}{2(m+1)}}  E \left[      \left( \sum_{l} \frac{|\zeta_{l}|^{p'}}{ (1+|l|^{2})^{sp'/2}}  \right)^{1/p'}  \right],~where~1/p'+1/p=1
\end{align*}
and since $\zeta$ is IID Gaussian and $sp'> 1$, $E \left[      \left( \sum_{l} \frac{|\zeta_{l}|^{p'}}{ (1+|l|^{2})^{sp'/2}}  \right)^{1/p'}  \right]$ is finite.  
	
This result and the Markov inequality implies the following uniform Law of Large Numbers result
	\begin{align}\label{eqn:Sobolev.Gaussian}
	\sup_{f \in \mathbb{W}^{s}_{p}(\mathbb{T})}  |n^{-1} \sum_{i=1}^{n} f(X_{i}) - E_{P}[f(X)] |  = O_{P}(n^{-\frac{1}{2} \frac{m}{m+1}}),
		\end{align}
	for slow $\theta$-mixing data. Thus generalizing  Glivenko-Cantelli results in Empirical processes theory (cf. \cite{VdV-W1996}) to ``slow" $\theta$-mixing processes. 
	
	However, the implied convergence rate can be slow when $m$ is close to one. Indeed, when $m=1$, the rate is $n^{-1/4}$, which is slower than the rate one would expect given that $\tau$ is ``almost summable". This feature shows a limitation of Proposition \ref{pro:bound.slow.gaussian} with respect to Corollary \ref{cor:bound.separation}: While the former is easier to implement, its implied bounds might be too conservative in some cases. The reason for this is that there is a trade-off between the norm defining the complexity measure, $L^{2r}(P)$, and the size of the scaling terms due to mixing; Proposition \ref{pro:bound.slow.gaussian} pushes this trade-off to one extreme by setting $r=1$. The next result leverages this observation and uses the bounds implied by Corollary \ref{cor:bound.separation} to shed more light on this issue and improve the rate,
	\begin{proposition}	
		Suppose $s>1$. Then, for any $r \in [1,+\infty]$, 
		\begin{align*}
			\sup_{f \in \mathbb{W}^{s}_{p}(\mathbb{T})}  |n^{-1} \sum_{i=1}^{n} f(X_{i}) - E_{P}[f(X)] |  = O_{P} \left( (\mathfrak{n}_{m,r'}(n))^{-1/2} \gamma_{2,1}( \mathbb{W}^{s}_{p}(\mathbb{T}) , ||.||_{L^{2r}(P)})   \right).
		\end{align*}
		with $r'$ such that $1/r'+1/r=1$ and where $n \mapsto \mathfrak{n}_{m,r'}(n) : = n^{ \frac{m(1+r')}{r'(1+m)} }$ if $r'>m$ and $n \mapsto \mathfrak{n}_{1,1}(n) : = n/(\log n) $ if $r'=m=1$.\footnote{If $r = \infty$, then we define $r'=1$ and vicerversa. If $r'=\infty$, then $\frac{m(1+r')}{r'(1+m)}$ is defined to be $\frac{m}{1+m}$.}
	\end{proposition}
	
\begin{proof}
	 Corollary \ref{cor:bound.separation} with $a=b=1$ and the Markov inequality imply
		\begin{align}\label{eqn:example.slow.1}
		\sup_{f \in \mathbb{W}^{s}_{p}(\mathbb{T})}  |n^{-1} \sum_{i=1}^{n} f(X_{i}) - E_{P}[f(X)] |  = O_{P} \left( \frac{ \bar{\tau}^{-1}(1/n)  }{n} \gamma_{1,1}( \mathbb{W}^{s}_{p}(\mathbb{T}) , ||.||_{L^{\infty}})    +    \frac{ \Theta_{r}(\bar{\tau}^{-1}(1/n))}{\sqrt{n} } \gamma_{2,1}( \mathbb{W}^{s}_{p}(\mathbb{T}) , ||.||_{L^{2r}(P)})   \right),
	\end{align}
	 for any $r \in [1,+\infty]$.

	 Observe that, for any $q \in \mathbb{N}$, $\int_{0}^{1} |\min\{ \alpha^{-1}(u), q\} |^{r'} du = \int_{0}^{1} |\min\{q,u^{-\frac{1}{m} }\}|^{r'} du =   q^{r'-m} + \int_{q^{-m}}^{1} u^{-r'/m} du$, and thus 
	$  \int_{0}^{1} |\min\{ \alpha^{-1}(u), q\} |^{r'} du =O (q^{r'-m} )$ for $r'>m$ and $ \int_{0}^{1} |\min\{ \alpha^{-1}(u), q\} |^{r'} du=O(\log q)$ for $r'=m$.  This result and the fact that $\bar{\tau}^{-1}(1/n) = O(  n^{1/(1+m)}   ) $ imply $\Theta_{r}(\bar{\tau}^{-1}(1/n)) = O(n^{\frac{r'-m}{2(1+m)r'}})$ for $r'>m$ and $\Theta_{r}(\bar{\tau}^{-1}(1/n))  = O(\sqrt{\log n})$ for $r'=m$. Also, $ \frac{ \bar{\tau}^{-1}(1/n)  }{\sqrt{n}} = O(n^{\frac{1-m}{2(1+m)}})$. 	
	Since $r' \geq 1$ these results imply that $ \frac{ \bar{\tau}^{-1}(1/n)  }{\sqrt{n}} = o\left( \Theta_{r}(\bar{\tau}^{-1}(1/n))  \right) $. 
	
	These results, the fact that $\gamma_{1,1}(\mathbb{W}^{s}_{p}(\mathbb{T}), ||.||_{L^{\infty}}) \leq \sum_{l} 2^{l-ls} < \infty$ (by expression \ref{eqn:gamma.entropy.bound} and  Remark \ref{rem:bound.gamma1} in Appendix \ref{app:bound.Dudley}), and expression \ref{eqn:example.slow.1} imply that
	\begin{align}\label{eqn:example.slow.2}
		\sup_{f \in \mathbb{W}^{s}_{p}(\mathbb{T})}  |n^{-1} \sum_{i=1}^{n} f(X_{i}) - E_{P}[f(X)] |  = O_{P} \left(  n^{-\frac{1}{2} }	\Theta_{r}(\bar{\tau}^{-1}(1/n))   \gamma_{2,1}( \mathbb{W}^{s}_{p}(\mathbb{T}) , ||.||_{L^{2r}(P)})   \right),
	\end{align}
	for any $r \in [1,+\infty]$. Which in turn implies 
			\begin{align*}
		\sup_{f \in \mathbb{W}^{s}_{p}(\mathbb{T})}  |n^{-1} \sum_{i=1}^{n} f(X_{i}) - E_{P}[f(X)] |  = O_{P} \left(  n^{-\frac{1}{2} \left( \frac{(1+r')m }{r'(1+m)} \right)  }  \gamma_{2,1}( \mathbb{W}^{s}_{p}(\mathbb{T}) , ||.||_{L^{2r}(P)})   \right),
	\end{align*}
	for any $r \in [1,+\infty]$ and $r'>m$; and 
	 			\begin{align*}
	 	\sup_{f \in \mathbb{W}^{s}_{p}(\mathbb{T})}  |n^{-1} \sum_{i=1}^{n} f(X_{i}) - E_{P}[f(X)] |  = O_{P} \left(  n^{-\frac{1}{2}  } \sqrt{\log n}  \gamma_{2,1}( \mathbb{W}^{s}_{p}(\mathbb{T}) , ||.||_{L^{2r}(P)})   \right),
	 \end{align*}
 for $r'=m=1$. 
	 
%
%
\end{proof}

	The rate from Proposition \ref{pro:bound.slow.gaussian} is recovered when $r=1$ (and $r'=\infty$). However, this rate can be improved by setting $r>1$. For instance, setting $r=\infty$ (and $r'=1$) implies, by expression \ref{eqn:gamma.entropy.bound} and Sobolev embedding results (e.g. \cite{EdmundsTriebel1996}), that $\gamma_{2,1}( \mathbb{W}^{s}_{p}(\mathbb{T}) , ||.||_{L^{\infty}}) < \infty$. So, the previous proposition yields a rate of $O_{P} \left( n^{-\frac{m}{1+m} }   \right)$ for $m<1$, which is faster than the one obtained in expression \ref{eqn:Sobolev.Gaussian}, especially when $m$ is close to one. On the other hand, if $m$ is close to zero --- or if in the application at hand, consistency regardless of the rate is enough --- then the previous proposition does not present a significant improvement over expression \ref{eqn:Sobolev.Gaussian}.	
	
	Finally, the quantity $\mathfrak{n}_{m,r'}(n)$ clearly illustrates the impact of the ``slow mixing" on the concentration rate. Indeed, $\mathfrak{n}_{m,r'}(n)$ and can be seen as an adjusted sample size, which modifies the actual sample size $n$ by incorporating the mixing structure, and generalizes the rate of $n/m$ obtained for the $m$-dependent case in expression \ref{eqn:bound.mdepen}. While this adjusted sample size diverges with the original one, it does so slower. $\triangle$ 
\end{example}

\section{Proofs}
\label{sec:proofs}	
	
\subsection{Proof of Theorem \ref{thm:equi}}
\label{sec:proof.main}

For the proof of this theorem it is useful to expand the notion of admissible partition sequence for $\mathcal{F}$ to a sequence of partitions $\mathcal{T}^{\infty}(c,d) : = (\mathcal{T}_{l}(c,d))_{l \in \mathbb{N}_{0}}$ which is increasing and $card \mathcal{T}_{0}(c,d) =1$  and  $card \mathcal{T}_{l}(c,d) \leq 2^{d 2^{l/c}}$ where $c > 0$ and $d >0$. The parameter $c$ is the order --- indeed when $d=1$, $\mathcal{T}_{l}(c) : = \mathcal{T}_{l}(c,1)$ coincides with the one in the definition \ref{def:MoC} ---, $d$ is an extra parameter that is convenient for establishing this proof. We call $(c,d)$ the order tuple of partition $\mathcal{T}^{\infty}(c,d)$.

For any order tuple $c,d>0$ and any admissible partition sequence for $\mathcal{F}$ given by $\mathcal{T}^{\infty}(c,d): = (\mathcal{T}_{k}(c,d))_{k \in \mathbb{N}_{0}}$, we construct a ``chain" from any $f$ to $f_{0}$ as follows. For any $k \in \mathbb{N}_{0}$, let $\pi_{k}f$ be an element  of $T(f,\mathcal{T}_{k}(c,d))$, where $T(f,\mathcal{T}_{k}(c,d))$ is the (unique) set in $\mathcal{T}_{k}(c,d)$ containing $f$. Setting $\pi_{0}f = f_{0}$, it follows that
\begin{align*}
	f-f_{0} = \sum_{l=1}^{\infty} \Delta_{l}f,~where~\Delta_{l}f : =  \pi_{l}f - \pi_{l-1}f.
\end{align*}

The following Lemma is the basis of the proof of Theorem \ref{thm:equi} and could be of independent interest as it gives a more primitive (albeit more cumbersome) bound than of the theorem; its proof is relegated to section \ref{sec:bound.general}. 

%

\begin{lemma}\label{lem:bound.general}
	For any $n \in \mathbb{N}$, any order tuple $c,d \geq 1$, any admissible partition sequence for $\mathcal{F}$, $\mathcal{T}^{\infty}(c,d)$, any real-valued sequence $(e(l))_{l \in \mathbb{N}_{0}}$, any $\mathbf{q} : \mathbb{N}_{0} \rightarrow \mathcal{Q}_{n}$, and any  $\ell : \mathbb{N}_{0} \rightarrow \mathbb{R}_{+}$  such that $\sup_{f \in \mathcal{F}} ||\Delta_{l} f ||_{L^{\infty}} \leq e(l)$ and $d 2^{l/c} \leq \ell(l)^{2} $ for all $l \in \mathbb{N}$, we obtain
	\begin{align*}
		\left \Vert  \sup_{f , f_{0}  \in \mathcal{F}} | G_{n}[f-f_{0}] |  \right \Vert_{L^{1}(P)} \leq  \mathbb{C}_{c,d}   \sup_{f \in \mathcal{F}}  \sum_{l=1}^{\infty} \sqrt{\frac{ \underline{n}(l) }{n}}  \left( ||\Delta_{l} f||_{2,\mathbf{q}(l)} \ell(l)   + ||\Delta_{l} f ||_{L^{\infty}}  \frac{\mathbf{q}(l)}{\sqrt{\underline{n}(l)}}  ( \ell(l)  )^{2}  +  \sqrt{\underline{n}(l) }  e(l) \tau (\mathbf{q}(l))  \right),
	\end{align*}
where $\mathbb{C}_{c,d} : = 2 + 48 \sum_{l=0}^{\infty} e^{-0.5 d 2 ^{l/c} }$ and $\underline{n}(l) : = J(n,\mathbf{q}(l)) \mathbf{q}(l)$ for any $l \in \mathbb{N}_{0}$ with $J(n,q)$ being the floor of $n/q$.
	
\end{lemma}

\begin{proof}[Proof of Theorem \ref{thm:equi}]
	
	To establish the desired result we employ Lemma \ref{lem:bound.general} with $\ell(l) = 2^{l/2} \sqrt{d}$ for all $l \in \mathbb{N}_{0}$, $d = 2^{1-1/c}$ and $c : = \min\{a,b\}$, and a $\mathcal{T}^{\infty}(c,d)$ and $(e(l))_{l \in \mathbb{N}_{0}}$ which we now construct.
	
	Let $\mathcal{A}^{\infty}(a)$ and $\mathcal{B}^{\infty}(b)$ be admissible partition sequences of order $a$ and $b$ respectively --- for now these are arbitrary partition sequences, we specify them later. For any $l \geq 1$, $\mathcal{T}_{l}(c,d)$ is comprised of sets of the form $A \cap B$ for $A \in \mathcal{A}_{l-1}(a)$ and $B \in \mathcal{B}_{l-1}(b)$. 
	
	Under this choice, for any $l \in \mathbb{N}$, $card \mathcal{T}_{l}(c,d) \leq 2^{2^{(l-1)/a} + 2^{(l-1)/b}}$. For $c = \min \{a,b\}$, $2^{(l-1)/a} + 2^{(l-1)/b} = 2^{l/c} 2^{-1/c} ( 2^{(l-1)(1/a-1/c)}  + 2^{(l-1)(1/b-1/c)}   ) \leq 2^{l/c} 2^{1-1/c} = d 2^{l/c}$. Hence, $card \mathcal{T}_{l}(c,d) \leq 2^{d2^{l/c}}$ for all $l \in \mathbb{N}_{0}$. 
	
	We now verify the condition $d 2^{l/c} \leq \ell(l)^{2}$  for all $l \in \mathbb{N}_{0}$. This condition is equivalent to $2^{l/c} \leq 2^{l}$ which is satisfied because $c  = \min\{a,b\} \geq 1$.

	Since $| \Delta_{l} f | \leq D(f,\mathcal{T}_{l-1})$ (as $\mathcal{T}_{l} \subseteq \mathcal{T}_{l-1}$) and $(\mathbf{q}_{n}(l))_{l} $ is non-increasing sequence, $\sum_{l=1}^{\infty}  2^{l/2} ||\Delta_{l}f||_{2,\mathbf{q}_{n}(l)} \leq 2 \sum_{l=0}^{\infty}  2^{l/2} || D(f,\mathcal{T}_{l}) ||_{2,\mathbf{q}_{n}(l)}  $ and similarly under $||.||_{\infty,\mathbf{q}_{n}(l)}$. Thus, we can choose $\mathcal{A}^{\infty}(a)$ and $\mathcal{B}^{\infty}(b)$ as the (approximately) optimal ones, i.e., such that  $\sup_{f \in  \mathcal{F}}   \sum_{l=1}^{\infty}  2^{l/2} ||\Delta_{l}f||_{2,\mathbf{q}_{n}(l)}  \leq 2 \gamma_{2,a}( \mathcal{F} , ||.||_{2,\mathbf{q}_{n}}  )(1+\epsilon)$ and $\sup_{f \in  \mathcal{F}}   \sum_{l=1}^{\infty}  2^{l} ||\Delta_{l}f||_{L^{\infty}} \mathbf{q}_{n}(l)  \leq 2 \gamma_{1,b}( \mathcal{F} ,  ||.||_{\infty,\mathbf{q}_{n}}   )(1+\epsilon)$ for some arbitrary small positive $\epsilon$. This readily implies that for all $f \in \mathcal{F}$ and  for all $l \in \mathbb{N}$ , $ ||\Delta_{l}f||_{L^{\infty}}  \leq 2 \frac{ \gamma_{1,b}( \mathcal{F} ,  ||.||_{\infty,\mathbf{q}_{n}}   ) }{2^{l}\mathbf{q}_{n}(l) } (1+\epsilon)$, and thus $e(l)$ can be taken as the RHS of the inequality.

	Therefore,  Lemma \ref{lem:bound.general} implies
	\begin{align*}
		\left \Vert  \sup_{f , f_{0}  \in \mathcal{F}} | G_{n}[f-f_{0}] |  \right \Vert_{L^{1}(P)} \leq  & \mathbb{C}_{c,d} \sup_{f \in \mathcal{F}}  \sum_{l=1}^{\infty}  \sqrt{\frac{ \underline{n}(l) }{n}}  \left( 2^{l/2} ||\Delta_{l} f||_{2,\mathbf{q}_{n}(l)}   + 2^{l} ||\Delta_{l} f ||_{L^{\infty}}  \frac{\mathbf{q}_{n}(l)}{\sqrt{\underline{n}(l)}}   \right) \\
		&  +  2.5 \mathbb{C}_{c,d}   \gamma_{1,b}( \mathcal{F} ,  ||.||_{\infty,\mathbf{q}_{n}}   ) \sum_{l=1}^{\infty} \sqrt{\underline{n}(l) }   \frac{  \tau (\mathbf{q}_{n}(l)) }{2^{l}\mathbf{q}_{n}(l) }  \\
		\leq &  \mathbb{C}_{c,d}  \sup_{f \in \mathcal{F}}  \sum_{l=1}^{\infty}  \left( 2^{l/2} ||\Delta_{l} f||_{2,\mathbf{q}_{n}(l)}   + 2^{l} ||\Delta_{l} f ||_{L^{\infty}}  \frac{\mathbf{q}_{n}(l)}{\sqrt{n}}   \right) \\
		&  +  2.5 \mathbb{C}_{c,d}   \gamma_{1,b}( \mathcal{F} ,  ||.||_{\infty,\mathbf{q}_{n}}   ) \sum_{l=1}^{\infty} \sqrt{ n }   \frac{  \tau (\mathbf{q}_{n}(l)) }{2^{l}\mathbf{q}_{n}(l) } \\
			\leq &  \mathbb{C}_{c,d}  \sup_{f \in \mathcal{F}}  \sum_{l=1}^{\infty}  \left( 2^{l/2} ||\Delta_{l} f||_{2,\mathbf{q}_{n}(l)}   + 2^{l} ||\Delta_{l} f ||_{L^{\infty}}  \frac{\mathbf{q}_{n}(l)}{\sqrt{n}}   \right) \\
		&  +  \mathbb{C}_{c,d}(3 + 2.5 \sqrt{2})   \frac{\gamma_{1,b}( \mathcal{F} ,  ||.||_{\infty,\mathbf{q}_{n}}   )}{\sqrt{n}},
	\end{align*}
where the second inequality holds because $ \underline{n}(l)\leq n$ for any $l \in \mathbb{N}_{0}$; the third inequality follows from the fact that, by expression \ref{eqn:qk}, $\mathbf{q}_{n}$ is such that $   \frac{  \tau (\mathbf{q}_{n}(l)) }{\mathbf{q}_{n}(l) } \leq 2^{l/2}/n$ for any $l \in \mathbb{N}_{0}$, which implies that $\sum_{l}  \frac{  \tau (\mathbf{q}_{n}(l)) }{2^{l}\mathbf{q}_{n}(l) } \leq  \sum_{l} \frac{2^{-l/2}}{n} \leq \mathbb{L}/n$. 

Since by our choices of partition, $\sup_{f \in  \mathcal{F}}   \sum_{l=1}^{\infty}  2^{l/2} ||\Delta_{l}f||_{2,\mathbf{q}_{n}(l)}  \leq 2 \gamma_{2,a}( \mathcal{F} , ||.||_{2,\mathbf{q}_{n}}  )(1+\epsilon)$  and  $\sup_{f \in  \mathcal{F}}   \sum_{l=1}^{\infty}  2^{l} ||\Delta_{l}f||_{L^{\infty}} \mathbf{q}_{n}(l)  \leq 2 \gamma_{1,b}( \mathcal{F} ,  ||.||_{\infty,\mathbf{q}_{n}}   )(1+\epsilon)$, the previous display readily implies 
	\begin{align*}
		\left \Vert  \sup_{f , f_{0}  \in \mathcal{F}} | G_{n}[f-f_{0}] |  \right \Vert_{L^{1}(P)} \leq \mathbb{L}_{a,b}  \left( \frac{\gamma_{1,b}( \mathcal{F} , ||.||_{\infty,\mathbf{q}_{n}}  ) }{\sqrt{n}}  + \gamma_{2,a}( \mathcal{F} , ||.||_{2,\mathbf{q}_{n}}  ) \right),
	\end{align*}
where $\mathbb{L}_{a,b} $ is a constant that depends only on $c,d$ --- and $c,d$ depend on $a,b$ --- and other universal constants.
\end{proof}

\subsubsection{Proof of Lemma \ref{lem:bound.general}}
\label{sec:bound.general}


For any $f \in \mathcal{F}$ and any $k \in \mathbb{N}_{0}$, let $\pi_{k}f$ be an element  of $T(f,\mathcal{T}_{k})$, where $T(f,\mathcal{T}_{k})$ is the (unique) set in $\mathcal{T}_{k}$ containing $f$.\footnote{To ease the notational burden, we leave the dependence of the partition on the order tuple $c,d$ implicit.} Let $\Pi_{k} \mathcal{F}$ be the subset of $\mathcal{F}$ of all such elements --- $\pi_{k}$ has the property that  $\pi_{k}f = \pi_{k}f'$ for all $f' \in T(f,\mathcal{T}_{k})$, thus $card \Pi_{k} \mathcal{F} \leq 2^{d 2^{k/c}}$. Henceforth, for any $l \in \mathbb{N}$, let $\Delta_{l} \mathcal{F} : = \{ g - g_{-1} \colon g \in \Pi_{l} \mathcal{F}~and~g_{-1} \in \Pi_{l-1} \mathcal{F}   \}$, and $B \Delta_{l} \mathcal{F} : = \{ g \in cone \Delta_{l} \mathcal{F} \colon ||g||_{L^{\infty}} \leq 1   \}$.

By setting $\pi_{0}f = f_{0}$, the following chain is constructed
\begin{align*}
	f-f_{0} = \sum_{l=1}^{\infty} (\pi_{l}f - \pi_{l-1}f),
\end{align*}
and so
\begin{align}\label{eqn:EP.chara.1}
 G_{n}[f-f_{0}] =  \sum_{l=1}^{\infty} G_{n}[\Delta_{l}f],~\forall f \in \mathcal{F}.
\end{align}

The goal is to couple the process $f \mapsto  G_{n}[\Delta_{l}f]$ with a block-independent one. However, since $n/q$ may not be an integer, we need to decompose the process as follows
\begin{align}\label{eqn:Gn.representation.1}
	G_{n}[\Delta_{l}f] = \sqrt{\frac{J(n,q)q}{n}} G_{J(n,q)q}[\Delta_{l}f]  + n^{-1/2} \sum_{j= J(n,q)q+1}^{n} \Delta_{l}f(X_{j})
\end{align}
where $J(n,q)$ is the floor of $n/q$, i.e., $J(n,q) : = \max\{  x \in \mathbb{N} \colon x \leq n/q \}$. This expression decomposes $G_{n}[f]$ into a ``reminder" part and a first part that can be ``coupled" with the block-independent process. The reminder term satisfies
\begin{align}\label{eqn:reminder.bound}
\left| 	n^{-1/2} \sum_{j= J(n,q)q+1}^{n} \Delta_{l}f(X_{j}) \right| \leq ||\Delta_{l}f||_{L^{\infty}} n^{-1/2} ( n - (J(n,q)q+1)  ) \leq  \frac{q}{\sqrt{n}} ||\Delta_{l}f||_{L^{\infty}}
\end{align}
 because  $0 \leq n-J(n,q)q \leq q$. 
 
For any $n,q,l \in \mathbb{N}$, the process $f \mapsto G_{J(n,q)q}[\Delta_{l}f] $ will be coupled with a block-independent process constructed using $(X^{\ast}_{i})_{i \in \mathbb{N}}$ given by\footnote{Throughout this section we use a double sub-index to denote the block-independent process --- $f \mapsto G^{\ast}_{n,q}[f]$ as opposed to just $f \mapsto G^{\ast}_{n}[f]$. This is to stress the dependence of the process on the length of the block $q$.}
\begin{align*}
 f \mapsto G^{\ast}_{J(n,q)q,q}[f] : = (J(n,q)q)^{-1/2} \sum_{j=0}^{J(n,q)-1} \left(   \delta f ( U^{\ast}_{j}(q) ) - E[ \delta f ( U^{\ast}_{j}(q) )] \right) \\
 with~\delta f ( U^{\ast}_{j}(q) )  : = \sum_{i=1}^{q} f(X^{\ast}_{qj+i}). 
\end{align*}
The existence and construction of process $(X^{\ast}_{i})_{i=-\infty}^{\infty}$ follow from known results and are relegated to Appendix \ref{app:coupling}. The process is such that the blocks $U^{\ast}_{j}(q) : =(X^{\ast}_{qj+1},...,X^{\ast}_{qj+q})$ have the same distribution as $U_{j}(q) : = (X_{qj+1},...,X_{qj+q})$ for each $j$, and $(U^{\ast}_{2j}(q))_{j=0}^{\infty}$ form an independent sequence and $(U^{\ast}_{2j+1}(q))_{j=0}^{\infty}$ form another independent sequence.

Henceforth, let $\underline{n}(l) : = J(n,\mathbf{q}(l)) \mathbf{q}(l)$ for any $l \in \mathbb{N}_{0}$. Then, by equations \ref{eqn:EP.chara.1} - \ref{eqn:reminder.bound}, and the triangle inequality, it follows that
 \begin{align}\notag
  \left \Vert \sup_{f , f_{0}  \in \mathcal{F}}	 G_{n}[f-f_{0}] \right \Vert_{L^{1}(P)} \leq  & \left \Vert \sup_{f \in \mathcal{F}}	\sum_{l=1}^{\infty} \sqrt{\frac{ \underline{n}(l) }{n}}  G^{\ast}_{\underline{n}(l),\mathbf{q}(l)}[\Delta_{l}f] \right \Vert_{L^{1}(P)} \\ \notag
  & +  \left \Vert \sup_{f \in \mathcal{F}}	\sum_{l=1}^{\infty} \sqrt{\frac{ \underline{n}(l) }{n}}  \left( G_{\underline{n}(l)}[\Delta_{l}f] -  G^{\ast}_{\underline{n}(l),\mathbf{q}(l)}[\Delta_{l}f] \right) \right \Vert_{L^{1}(\mathbb{P})}\\ \label{eqn:L1bound.master.1}
  & + \sup_{f \in \mathcal{F} }	\sum_{l=1}^{\infty} \frac{\mathbf{q}(l)}{\sqrt{n}} ||\Delta_{l}f||_{L^{\infty}}.
 \end{align}
($\mathbb{P}$ denotes the joint probability measure over the original process and the block-independent one.)
 
We now bound the first term in the RHS in expression \ref{eqn:L1bound.master.1},  $	\sup_{f \in \mathcal{F}} | \sum_{l=1}^{\infty}  \sqrt{\frac{ \underline{n}(l) }{n}} G^{\ast}_{\underline{n}(l),\mathbf{q}(l)}[\Delta_{l} f] | $. To do this, we employ an extension of the generic chaining approach proposed by Talagrand (cf. see \cite{talagrand2014}). One condition needed for this approach (cf. \cite{talagrand2014} expression 1.4) is a Bernstein-type inequality for $ G^{\ast}_{\underline{n}(l),\mathbf{q}(l)}[\Delta_{l} f]$ which is verified in the following lemma
\begin{lemma}\label{lem:bere}
	For any $g \in \mathcal{F}$ and any $n \in \mathbb{N}$ any $q \in \mathcal{Q}_{n}$ such that $n/q \in \mathbb{N}$ it follows that 
	\begin{align*}
		P  \left(	G^{\ast}_{n,q}[g] \geq 6 ||g||_{2,q} u  + \frac{16}{3} \frac{q}{\sqrt{n}} ||g||_{L^{\infty}} u^{2} \right) \leq 2 e^{-u^{2}},
	\end{align*}			
	for any $u > 0$.
\end{lemma}

\begin{proof}
	See Appendix \ref{app:main.proof}.
\end{proof}

We are now in a position to establish an exponential tail inequality for $\sup_{f \in \mathcal{F}} \left| \sum_{l=1}^{\infty}  \sqrt{\frac{ \underline{n}(l) }{n}} G^{\ast}_{\underline{n}(l),\mathbf{q}(l)}[\Delta_{l} f] \right|$.

\begin{lemma}\label{lem:bound.pr.talagrand}
	For any $n \in \mathbb{N}$, any order tuple $c,d \geq 1$, any admissible partition sequence for $\mathcal{F}$, $\mathcal{T}^{\infty}(c,d)$, any $\mathbf{q} : \mathbb{N}_{0} \rightarrow \mathcal{Q}_{n}$, and any  $\ell : \mathbb{N}_{0} \rightarrow \mathbb{R}_{+}$  such that $d 2^{l/c} \leq \ell(l)^{2} $ for all $l \in \mathbb{N}$,
	\begin{align*}
		&P \left(  \sup_{f \in \mathcal{F}}  \sum_{l=1}^{\infty} \sqrt{\frac{ \underline{n}(l) }{n}} |G^{\ast}_{\underline{n}(l),\mathbf{q}(l)}[\Delta_{l} f] | \leq v   \sup_{f \in \mathcal{F}}  \sum_{l=1}^{\infty}   \sqrt{\frac{\underline{n}(l)}{n}} \left( 6 ||\Delta_{l} f||_{2,\mathbf{q}(l)} \ell(l)   + \frac{16}{3} \frac{\mathbf{q}(l)}{\sqrt{ \underline{n}(l) }} ||\Delta_{l} f||_{L^{\infty}} (\ell(l)  )^{2}    \right)   \right) \\
		& \geq 1 - \left( 2 \sum_{l=0}^{\infty} 2^{  2d 2^{l/c}  }  e^{- 2\ell(l)^{2} } \right) e^{-0.5 v}
	\end{align*}
	for any $v \geq 4$.
\end{lemma}

\begin{proof}
		See Appendix \ref{app:main.proof}.
\end{proof}

From this previous lemma we can derive an expectation bound:

\begin{lemma}\label{lem:bound.E.talagrand}
	For any $n \in \mathbb{N}$, any order tuple $c,d \geq 1$, any admissible partition sequence for $\mathcal{F}$, $\mathcal{T}^{\infty}(c,d)$, any $\mathbf{q} : \mathbb{N}_{0} \rightarrow \mathcal{Q}_{n}$, and any  $\ell : \mathbb{N}_{0} \rightarrow \mathbb{R}_{+}$  such that $d 2^{l/c} \leq \ell(l)^{2} $ for all $l \in \mathbb{N}$,
	\begin{align*}
		 E_{P}  \left[    \sup_{f \in \mathcal{F}} \sum_{l=1}^{\infty}  \sqrt{\frac{ \underline{n}(l) }{n}}  |  G^{\ast}_{\underline{n}(l),\mathbf{q}(l)}[\Delta_{l} f] |    \right]  \leq \mathbb{M}  \sup_{f \in \mathcal{F}}  \sum_{l=1}^{\infty} \sqrt{\frac{ \underline{n}(l) }{n}}  \left( ||\Delta_{l} f||_{2,\mathbf{q}(l)} \ell(l)   + ||\Delta_{l} f ||_{L^{\infty}}  \frac{\mathbf{q}(l)}{\sqrt{\underline{n}(l)}}  ( \ell(l)  )^{2}   \right)  
	\end{align*}
	where $\mathbb{M} : = 12  \left(  \sum_{l=0}^{\infty} 2^{  2d 2^{l/c}  }  e^{- 2\ell(l)^{2} }  \right) \int_{0}^{\infty} e^{-0.5t } dt$.
\end{lemma}

\begin{proof}
		See Appendix \ref{app:main.proof}.
\end{proof}

By plugging in the bound in Lemma \ref{lem:bound.E.talagrand} in expression \ref{eqn:L1bound.master.1} one obtains
 \begin{align*}
	\left \Vert \sup_{f , f_{0}  \in \mathcal{F}}	 G_{n}[f-f_{0}] \right \Vert_{L^{1}(P)} \leq  & \mathbb{M}  \sup_{f \in \mathcal{F}}  \sum_{l=1}^{\infty} \sqrt{\frac{ \underline{n}(l) }{n}}  \left( ||\Delta_{l} f||_{2,\mathbf{q}(l)} \ell(l)   + ||\Delta_{l} f ||_{L^{\infty}}  \frac{\mathbf{q}(l)}{\sqrt{\underline{n}(l)}}  ( \ell(l)  )^{2}   \right)  \\ 
	& +  \left \Vert \sup_{f \in \mathcal{F}}	\sum_{l=1}^{\infty} \sqrt{\frac{ \underline{n}(l) }{n}}  \left( G_{\underline{n}(l)}[\Delta_{l}f] -  G^{\ast}_{\underline{n}(l),\mathbf{q}(l)}[\Delta_{l}f] \right) \right \Vert_{L^{1}(\mathbb{P})}\\ 
	& + \sup_{f \in \mathcal{F} }	\sum_{l=1}^{\infty} \frac{\mathbf{q}(l)}{\sqrt{n}} ||\Delta_{l}f||_{L^{\infty}}.
\end{align*}

Since $(\ell(l))^{2} \geq 1$, the last term in the RHS --- the ``reminder" term of the coupling --- is bounded by $ \sup_{f \in \mathcal{F} }	\sum_{l=1}^{\infty} \frac{\mathbf{q}(l)}{\sqrt{n}} ||\Delta_{l}f||_{L^{\infty}} \leq  \sup_{f \in \mathcal{F} }	\sum_{l=1}^{\infty} ||\Delta_{l}f||_{L^{\infty}} \frac{\mathbf{q}(l)}{\sqrt{n}} (\ell(l))^{2} $, which is dominated by the first term in the RHS. Thus,
 \begin{align}\notag
	\left \Vert \sup_{f , f_{0}  \in \mathcal{F}}	 G_{n}[f-f_{0}] \right \Vert_{L^{1}(P)} \leq  & (1 + \mathbb{M})  \sup_{f \in \mathcal{F}}  \sum_{l=1}^{\infty} \sqrt{\frac{ \underline{n}(l) }{n}}  \left( ||\Delta_{l} f||_{2,\mathbf{q}(l)} \ell(l)   + ||\Delta_{l} f ||_{L^{\infty}}  \frac{\mathbf{q}(l)}{\sqrt{\underline{n}(l)}}  ( \ell(l)  )^{2}   \right)  \\
	& +  \left \Vert \sup_{f \in \mathcal{F}}	\sum_{l=1}^{\infty} \sqrt{\frac{ \underline{n}(l) }{n}}  \left( G_{\underline{n}(l)}[\Delta_{l}f] -  G^{\ast}_{\underline{n}(l),\mathbf{q}(l)}[\Delta_{l}f] \right) \right \Vert_{L^{1}(\mathbb{P})}. \label{eqn:L1bound.master-2}
\end{align}

We now bound the second term in the RHS. By the condition in the lemma, $||\Delta_{l}f||_{L^{\infty}} \leq e(l)$ for any $f \in \mathcal{F}$ and any $l \in \mathbb{N}$. So, 
\begin{align*}
		 \left \Vert \sup_{f \in \mathcal{F}}    \left|  G_{\underline{n}(l)}[\Delta_{l}f] -  G^{\ast}_{\underline{n}(l),\mathbf{q}(l)}[\Delta_{l}f] \right| \right \Vert_{L^{1}(\mathbb{P})} 	= &	 \left \Vert \sup_{h \in \Delta_{l} \mathcal{F}}    \left|  G_{\underline{n}(l)}[h] -  G^{\ast}_{\underline{n}(l),\mathbf{q}(l)}[h] \right| \right \Vert_{L^{1}(\mathbb{P})} 	\\
		 \leq &	e(l) \left \Vert \sup_{g \in B \Delta_{l} \mathcal{F}}   \left|  G_{\underline{n}(l)}[g] -  G^{\ast}_{\underline{n}(l),\mathbf{q}(l)}[g] \right| \right \Vert_{L^{1}(\mathbb{P})} ,
\end{align*}
where $B \Delta_{l} \mathcal{F} : = \{ f \in cone  \Delta_{l} \mathcal{F}  \colon ||f||_{\infty} \leq 1    \}$.
By Lemma \ref{lem:tau-bdd} --- in which $\underline{n}(l)$ plays the role of $n$, $\mathcal{B} = B \Delta_{l} \mathcal{F}$, and $q = \mathbf{q}(l)$ --- the norm in the last term can be bounded by $\sqrt{\underline{n}(l) }   \tau_{B \Delta _{l} \mathcal{F}} (\mathbf{q}(l))$, thereby obtaining
\begin{align}
	 \left \Vert \sup_{f \in \mathcal{F}}    \left|  G_{\underline{n}(l)}[\Delta_{l}f] -  G^{\ast}_{\underline{n}(l),\mathbf{q}(l)}[\Delta_{l}f] \right| \right \Vert_{L^{1}(\mathbb{P})} 	
	 \leq &  e(l) \sqrt{\underline{n}(l) }  \tau_{B \Delta _{l} \mathcal{F}} (\mathbf{q}(l)). \label{eqn:L1bound.coupling.1}
\end{align}

This result and expression \ref{eqn:L1bound.master-2} imply that, for any $n \in \mathbb{N}$, any order tuple $c,d \geq 1$, any admissible partition sequence for $\mathcal{F}$, $\mathcal{T}^{\infty}(c,d)$, any real-valued sequence $(e(l))_{l \in \mathbb{N}_{0}}$, any $\mathbf{q} : \mathbb{N}_{0} \rightarrow \mathcal{Q}_{n}$, and any  $\ell : \mathbb{N}_{0} \rightarrow \mathbb{R}_{+}$  such that $\sup_{f \in \mathcal{F}} ||\Delta_{l} f ||_{L^{\infty}} \leq e(l)$ and $d 2^{l/c} \leq \ell(l)^{2} $ for all $l \in \mathbb{N}$, it follows that
\begin{align*}
	\left \Vert  \sup_{f , f_{0}  \in \mathcal{F}} | G_{n}[f-f_{0}] |  \right \Vert_{L^{1}(P)} \leq \mathbb{M}_{2}  \sup_{f \in \mathcal{F}}  \sum_{l=1}^{\infty} \sqrt{\frac{ \underline{n}(l) }{n}}  \left( ||\Delta_{l} f||_{2,\mathbf{q}(l)} \ell(l)   + ||\Delta_{l} f ||_{L^{\infty}}  \frac{\mathbf{q}(l)}{\sqrt{\underline{n}(l)}}  ( \ell(l)  )^{2}  +  \sqrt{\underline{n}(l) }  e(l)  \tau_{B \Delta _{l} \mathcal{F}} (\mathbf{q}(l))  \right),
\end{align*}
with $\mathbb{M}_{2} : = 1 + 12  \left( 2 \sum_{l=0}^{\infty} 2^{  2d 2^{l/c}  }  e^{- 2\ell(l)^{2} }  \right) \geq 1 + \mathbb{M}$.

We now establish that $  \tau_{B \Delta_{l} \mathcal{F}}(.) \leq 2  \tau_{B \mathcal{F}}(.) = 2 \tau(.)$ for any $l \in \mathbb{N}$. The last equality follows from the definition of $\tau$ in expression \ref{eqn:defn.tautheta} --- recall that for any set $S$, $BS : = \{ f \in cone S \colon ||f||_{L^{\infty}} \leq 1  \}$. Hence, we only need to show the inequality. To do this, we show that $  \tau_{B \Delta_{l} \mathcal{F}}(.) \leq 2 \tau_{B \Pi_{l} \mathcal{F} }(.) \leq 2\tau_{B \mathcal{F}}(.) $ where, recall, $\Pi_{l} \mathcal{F} : = \{ \pi_{l}f \colon f \in \mathcal{F} \}$. Since $B \mapsto \tau_{B}$ is non-decreasing (with respect to the set inclusion partial ordering; see Lemma \ref{lem:tau.properties}) and $\Pi_{l} \mathcal{F} \subseteq \mathcal{F}$, the last inequality follows. To show the first one, suppose, $\Lambda(\Delta_{l} \mathcal{F}) \subseteq \Lambda(2 \Pi_{l}\mathcal{F})$, then, by the definition of $\tau$, $\tau_{\Delta_{l} \mathcal{F}} \leq \tau_{2 \Pi_{l} \mathcal{F}} \leq 2\tau_{\Pi_{l} \mathcal{F}}$, so the desired inequality follows.  We now show $\Lambda(\Delta_{l} \mathcal{F}) \subseteq \Lambda(2 \Pi_{l}\mathcal{F})$. Take any $g \in \Lambda(\Delta_{l} \mathcal{F})$, it follows that $|g(x)-g(y)| \leq \sup_{h \in \Delta_{l} \mathcal{F}} |h(x)-h(y)|$. Since $h$ has to be of the form $f_{l}-f_{-l}$ for $f_{l},f_{-l}$ in $\Pi_{l}\mathcal{F}$, it follows that  $|g(x)-g(y)| \leq \sup_{h \in 2 \Pi_{l}\mathcal{F}} |h(x)-h(y)|$. Hence $g \in \Lambda(2 \Pi_{l}\mathcal{F})$.

Therefore, 
\begin{align*}
	\left \Vert  \sup_{f , f_{0}  \in \mathcal{F}} | G_{n}[f-f_{0}] |  \right \Vert_{L^{1}(P)} \leq 2 \mathbb{M}_{2}  \sup_{f \in \mathcal{F}}  \sum_{l=1}^{\infty} \sqrt{\frac{ \underline{n}(l) }{n}}  \left( ||\Delta_{l} f||_{2,\mathbf{q}(l)} \ell(l)   + ||\Delta_{l} f ||_{L^{\infty}}  \frac{\mathbf{q}(l)}{\sqrt{\underline{n}(l)}}  ( \ell(l)  )^{2}  +  \sqrt{\underline{n}(l) }  e(l)  \tau (\mathbf{q}(l))  \right).
\end{align*}
Thus, the desired result follows because $\sum_{l=0}^{\infty} 2^{  2d 2^{l/c}  }  e^{- 2\ell(l)^{2} }  \leq \sum_{l=0}^{\infty} e^{- m d 2^{l/c}}$ where $m : = 2(1-\ln 2) \geq 1/2$.

\subsection{Proof of Corollary \ref{cor:bound.separation}}
\label{sec:bound.separation}

\begin{proof}[Proof of Corollary \ref{cor:bound.separation}]

	Throughout the proof, let $r'$ be such that $1/r'+1/r=1$ with the convention that $r' = + \infty$ when $r=1$ and vicerversa. 
	
	We first show that 
	\begin{align}\label{eqn:bound.sep.proof.1}
		\left\Vert \sup_{f , f_{0}  \in \mathcal{F}}	 G_{n}[f-f_{0}]   \right\Vert_{L^{1}(P)} \leq  \mathbb{L}_{a,b} \frac{q_{n}}{\sqrt{n}}   \gamma_{1,b}( \mathcal{F} , ||.||_{L^{\infty}}  )  +  \gamma_{2,a}( \mathcal{F} , ||.||_{2,q_{n}}  ) 
	\end{align}
	for an integer $q_{n} : = \min\{ s \leq n \colon  \tau(s)n  \leq s  \}$. 
	
	There are two ways of establishing this result. One is to invoke Theorem \ref{thm:equi} and majorized the RHS using the fact that $\mathbf{q}_{n}(0) \geq \mathbf{q}_{n}(l)$ for all $l \in \mathbb{N}_{0}$ --- it is clear that $\mathbf{q}_{n}(0) = q_{n}$. Another way, which is perhaps more transparent and better illustrates the ideas behind the proof, is to show this result by first coupling and then applying a chaining arguments to the coupled process only. 
	
	That is, by triangle inequality, for any $q \leq n$,
	 \begin{align}\notag
		\left \Vert \sup_{f , f_{0}  \in \mathcal{F}}	 G_{n}[f-f_{0}] \right \Vert_{L^{1}(P)} \leq  & \sqrt{\frac{ J(n,q)q }{n}}   \left \Vert \sup_{f,f_{0} \in \mathcal{F}}  G^{\ast}_{J(n,q)q,q}[f-f_{0}] \right \Vert_{L^{1}(P)} \\ \notag
		& +  \sqrt{\frac{ J(n,q)q }{n}}  \left \Vert \sup_{f,f_{0} \in \mathcal{F}}  \left( G_{J(n,q)q}[f-f_{0}] -  G^{\ast}_{J(n,q)q,q}[f-f_{0}] \right) \right \Vert_{L^{1}(\mathbb{P})}\\ \label{eqn:L1bound.q0}
		& + \sup_{f,f_{0} \in \mathcal{F} }\frac{q}{\sqrt{n}} ||f-f_{0}||_{L^{\infty}},
	\end{align}
where the last term stems from the reminder term when $n/q$ is not an integer. 
	
	By applying the chaining arguments presented in Lemmas \ref{lem:bound.pr.talagrand} and \ref{lem:bound.E.talagrand} but with $\mathbf{q}(.) = q$ to the process $G^{\ast}_{J(n,q)q,q}$, the first term in the RHS is majorized (up to constants) by $$ \sup_{f \in \mathcal{F}}  \left\{   \sqrt{\frac{ J(n,q)q }{n}}   \sum_{l=0}^{\infty}  2^{l/2}    || D(f,\mathcal{T}_{l})  ||_{2,q}  + \frac{q}{\sqrt{n}}  \sum_{l=0}^{\infty} 2^{l}  ||D(f,\mathcal{T}_{l})  ||_{L^{\infty}}   \right\} .  $$
	
	By following the same steps that lead to inequality \ref{eqn:L1bound.coupling.1}, the second term in the RHS in expression \ref{eqn:L1bound.q0} is majorized (up to constants) by  $\frac{ J(n,q)q } {\sqrt{n}} \sup_{f,f_{0} \in \mathcal{F} }  ||f-f_{0}||_{L^{\infty}}  \tau(q)$. 
	
	These bounds and the facts that $ J(n,q)q \leq n$ and $\sup_{f,f_{0} \in \mathcal{F} }  ||f-f_{0}||_{L^{\infty}} = \sup_{f \in \mathcal{F} }  ||D(f,\mathcal{T}_{0})||_{L^{\infty}} \leq  \sum_{l=0}^{\infty} 2^{l}  ||D(f,\mathcal{T}_{l})  ||_{L^{\infty}} $, imply that 
	 \begin{align}\notag
	\left \Vert \sup_{f , f_{0}  \in \mathcal{F}}	 G_{n}[f-f_{0}] \right \Vert_{L^{1}(P)} \leq  \mathbb{L}_{a,b} \sup_{f \in \mathcal{F}}  \left\{     \sum_{l=0}^{\infty}  2^{l/2}    || D(f,\mathcal{T}_{l})  ||_{2,q}  + \left( \frac{q}{\sqrt{n}} + \sqrt{n} \tau(q) \right) \sum_{l=0}^{\infty} 2^{l}  ||D(f,\mathcal{T}_{l})  ||_{L^{\infty}}   \right\}.
\end{align}	

The choice $q=q_{n}$ precisely yields  $ \frac{q}{\sqrt{n}} + \sqrt{n} \tau(q) \leq  \mathbb{L} \frac{q}{\sqrt{n}}$ and thus, by optimizing over the partition sequence, the desired result follows. 	Expression \ref{eqn:bound.sep.proof.1} is thus proven. 

We now majorize $||.||_{2,q_{n}} $. To do this, by H\"{o}lder inequality, 
	\begin{align*}
	||f||_{2,q} \leq \sqrt{||\mu_{q}||_{L^{r'}([0,1],Leb)}  ||Q^{2}_{f}||_{L^{r}([0,1],Leb)}} = \sqrt{||\mu_{q}||_{L^{r'}([0,1],Leb)} } ||f||_{L^{2r}(P)} 
	\end{align*}
	for any $f \in \mathcal{F}$, any $q \in \mathbb{N}$, and any $r \in [1,+\infty]$ (with the convention that $L^{2r}(P)=L^{\infty}(P)$).

	We now show that for $r' < \infty$, $\sqrt{ || \mu_{q} ||_{L^{r'}([0,1],Leb)} }   \leq  \sqrt{  \left(  \int_{0}^{1} | \min\{ \alpha^{-1}(u), q  \} |^{r'} du   \right)^{1/r'} + 1 }$ and for $r' = \infty$, $\sqrt{  || \mu_{q} ||_{L^{r'}([0,1],Leb)} }  \leq \sqrt{1+q}$,  for any $q \in \mathbb{N}$ and thus establish the desired result. The last inequality follows directly from the construction of $\mu_{q}$ in expression \ref{eqn:muq}. This first inequality follows  by Lemma  \ref{lem:mu.bound} in Appendix \ref{app:supp.lemmas}. So
	\begin{align*}
	|| \mu_{q} ||_{L^{r'}([0,1],Leb)} \leq  \left( \int_{0}^{1} | \min\{ \alpha^{-1}(u), q  \}  +1 |^{r'} du  \right)^{1/r'}    \leq & \left(  \int_{0}^{1} | \min\{ \alpha^{-1}(u), q  \} |^{r'} du   \right)^{1/r'} + 1.
	\end{align*}
%
%
%
%
\end{proof}

\small
\bibliographystyle{plain}
\bibliography{myref}

\appendix 

\pagebreak
\begin{center}
	{\huge{Appendix}}
\end{center}
\bigskip

\section{Coupling Technique} 
\label{app:coupling}

We first formally construct the sequence $(X^{\ast}_{i})_{i=-\infty}^{\infty}$ using the results in \cite{dedecker2006inequalities,dedecker2004coupling,dedecker2005new}; we do this here for completeness, there are no conceptual innovations in this section. We then use these results to provide a coupling bound for our empirical process. 

Throughout, we use $\sigma(X)$ to denote the $\sigma$-algebra generated by a random variable $X$. 


\paragraph{$\tau$-coupling results (\cite{dedecker2006inequalities,dedecker2004coupling,dedecker2005new}).}Let $\mathsf X$ be a Polish space. Let $d$ be a distance on $\mathsf X$
(the space $\mathsf X$ need not be Polish with respect to $d$).
Let $\Lambda_1(\mathsf X;d)$ be the set of $1$--Lipschitz functions from $\mathsf X$ to $\mathbb R$
with respect to $d$. Assume that $d$ satisfies the Kantorovich--Rubinstein duality condition
\begin{align}
	\label{eq:KR}
	d(x,y)
	=
	\sup_{f\in\Lambda_1(\mathsf X;d)} |f(x)-f(y)|.
\end{align}

Let $(\Omega,\mathcal A,\mathbb P)$ be a probability space.
We say that an $\mathsf X$--valued random variable $X$ belongs to $L^r(\mathsf X;d)$ for some $r \in [1,\infty)$
if $d(\cdot,s)\in L^r(\mathbb R)$ for some (and therefore any) $s\in\mathsf X$.

For any $X\in L^1(\mathsf X;d)$ and any sub-$\sigma$--algebra $\mathcal M\subset\mathcal A$,
let $P_{X\mid\mathcal M}$ be a conditional distribution of $X$ given $\mathcal M$
and let $P_X$ be the (unconditional) distribution of $X$.
Define the dependence coefficient
\begin{align}
	\label{eq:tau-def}
	\tau_{d}(\mathcal M;X)
	:=
	\left\|
	\sup_{f\in\Lambda_1(\mathsf X;d)}
	\left|
	\int f(x)\,P_{X\mid\mathcal M}(dx)
	-
	\int f(x)\,P_X(dx)
	\right|
	\right\|_{L^{1}(\mathbb{P})}
\end{align}

The next lemma is Lemma 1 in \cite{dedecker2006inequalities}.

\begin{lemma}[Coupling representation of $\tau_{d}$]
	\label{lem:coupling}
	Assume that there exists a random variable $\zeta$ uniformly distributed over $[0,1]$,
	independent of $\sigma(X)\vee \mathcal M$.
	Then there exists a random variable $X^\ast$ measurable with respect to
	\begin{align}
		\sigma(\zeta)\vee\sigma(X)\vee\mathcal M,
	\end{align}
	distributed as $X$ and independent of $\mathcal M$, such that,
	\begin{align}
		\label{eq:coupling-eq}
		\tau_d(\mathcal M;X)
		&=
		\mathbb E\big[d(X,X^\ast)\big].
	\end{align}
\end{lemma}

\paragraph{Construction of $(X_i^\ast)_{i=-\infty}^{+\infty}$ from Lemma \ref{lem:coupling}.} Throughout the rest of this section, fix any $n \ge q$ in $\mathbb{N}$ such that $n/q \in \mathbb{N}$. Let 
\begin{align*}
	U_{j}(q) : = (X_{(j-1)q+1},...,X_{(j-1)q+q})
\end{align*}
for any $j \in [n/q]$.\footnote{For an integer $n$, $[n] :  = \{1,...,n\} $.} Observe that $U_{j}(q) \in \mathbb{X}^{q}$ for  any $j \in [n/q]$. Let
\begin{align}
	\mathcal U_{j-2}(q) := \sigma(U_1(q),\ldots,U_{j-2}(q)), \qquad j\ge 3,
\end{align}
and $(\zeta_j)_{j\ge 3}$ be IID random variables, each uniformly distributed on $[0,1]$,
and independent of $\sigma(U_1(q),U_2(q),\ldots)$. Henceforth, we omit the dependence on $q$ on $U_{j}$, $\mathcal U_{j-2}$, and other analogous quantities. 

Fix $j \in [n/q]$ and, following the notation in the section above, let $\mathsf X = \mathbb{X}^{q}$, $X = U_{j}$, $\mathcal{M} = \mathcal{U}_{j-2}$, $\zeta = \zeta_{j}$, and $d_{q}$ be the notion of distance over $\mathbb X^{q}$. By Lemma \ref{lem:coupling} we obtain an $\mathbb X^{q} $--valued random variable $U_j^\ast = : (X^{\ast}_{(j-1)q+1},...,X^{\ast}_{(j-1)q+q})$ such that:
\begin{align}
	& U_j^\ast \stackrel{d}{=} U_j,~U_j^\ast \perp\!\!\!\perp \mathcal U_{j-2},~U_j^\ast \in \sigma(\zeta_j)\vee\sigma(U_j)\vee\mathcal U_{j-2},~and \\ \label{eqn:coupling.Hq}
	& \mathbb E\big[ d_{q}(U_j,U_j^\ast) \big] = \tau_{d_{q}}(\mathcal U_{j-2};U_j).
\end{align}

The above ``two-step'' decoupling implies that the even subsequence $(U_{2k}^\ast)_{k\ge 1}$
is mutually independent (and likewise for the odd subsequence $(U_{2k-1}^\ast)_{k\ge 1}$).
Indeed, fix $k<l$. Since $2k\le 2l-2$, we have
\begin{align}
	U_{2k}^\ast \in \sigma(U_1,\ldots,U_{2l-2})=\mathcal U_{2l-2}.
\end{align}
But by construction
\begin{align}
	U_{2l}^\ast \perp\!\!\!\perp \mathcal U_{2l-2},
\end{align}
hence $U_{2l}^\ast \perp\!\!\!\perp U_{2k}^\ast$. Since this holds for every earlier even index,
the family $(U_{2k}^\ast)_{k\ge 1}$ is independent. The argument for odd indices is identical.

The coupling used in this paper is constructed using metrics of the form
\begin{align}\label{eqn:metric.F.q}
	(x_{1},...,x_{q} ,x'_{1},...,x'_{q}) \mapsto d_{q,\mathcal B}(x_{1},...,x_{q} ,x'_{1},...,x'_{q} ) : = \sum_{m=1}^{q} \sup_{f \in \mathcal{B} } | f(x_{m}) - f(x'_{m})|
\end{align}
for $\mathcal{B} \subseteq cone \mathcal{F}$ provided $\mathbb E [\sup_{f \in \mathcal{B}} |f(X_{0})| ] < \infty$.  Under this (pseudo) metric, the  expression \ref{eqn:coupling.Hq} yields 
\begin{align}\label{eqn:coupling.Hq.1}
	\mathbb E\left[ \sum_{m=1}^{q} \sup_{f \in \mathcal{B}} | f(X_{(j-1)q+m}) - f(X^{\ast}_{(j-1)q+m})|  \right] = \tau_{d_{q,\mathcal B}}(\sigma(U_1,\ldots,U_{j-2});U_j).
\end{align}

\subsection{Coupling results} 
\label{app:coupling.bounds}

The next results apply the construction in the previous section to the generate a coupling bound for the empirical process $G_{n}$ over some $\mathcal{B} \subseteq cone \mathcal{F}$.  To do this, consider $d_{q,\mathcal{B}}$ as in expression \ref{eqn:metric.F.q}. Given this and any $n \in \mathbb{N}$ and $q \in \mathcal{Q}_{n}$ such that $n/q \in \mathbb{N}$, let $(X^{\ast}_{i})_{i \in \mathbb{N}}$ be constructed as above and 
\begin{align*}
	G^{\ast}_{n,q}[f] = n^{-1/2} \sum_{j=0}^{n/q-1} \sum_{l=1}^{q}  \{ f(X^{\ast}_{qj+l})  - E_{P}[f(X^{\ast})] \}.  
\end{align*}

Throughout, for integers $a\le b$, let
\[
\mathcal M_a^b := \sigma(X_i:\ a\le i\le b),\qquad
\mathcal M_{-\infty}^b := \mathcal M^{b} : = \sigma(X_i:\ i\le b),\qquad
\mathcal M_a^{+\infty}:=  \mathcal M_{a} : =  \sigma(X_i:\ i\ge a).
\]

\begin{lemma}
	\label{lem:tau-bdd.0}
	Let $\mathcal{B} \subseteq cone \mathcal{F}$. For any $n \in \mathbb{N}$ and $q \in \mathcal{Q}_{n}$ such that $n/q \in \mathbb{N}$, 
	\begin{align*}
		\left \Vert \sup_{f \in \mathcal{B} } | G_{n}[f]  - G^{\ast}_{n,q}[f]     |  \right \Vert_{L^{1}(\mathbb{P})} \leq  n^{-1/2} \sum_{j=0}^{n/q-1}   \tau_{d_{q,\mathcal B}}(\mathcal{M}_{1}^{jq-q};(X_{jq+1},\ldots ,X_{jq+q} )).
	\end{align*}
	
\end{lemma}

\begin{proof}[Proof of Lemma \ref{lem:tau-bdd.0}]
	
	For any $n \in \mathbb{N}$ and $q \leq n$ such that $n/q$ is an integer
	\begin{align*}
		G_{n}[f]  = n^{-1/2} \sum_{j=0}^{n/q-1} \sum_{l=1}^{q}  \{ f(X_{qj+l}) -   E_{P}[f(X)]   \}  .
	\end{align*}
	
	Hence, 
	\begin{align*}
		G_{n}[f]  - G^{\ast}_{n,q}[f]  = n^{-1/2} \sum_{j=0}^{n/q-1} \sum_{l=1}^{q} \left(  f(X_{qj+l}) - f(X^{\ast}_{qj+l})     \right) .
	\end{align*}
	Then, by definition of $d_{q,\mathcal{B}}$,
	\begin{align*}
		E_{\mathbb{P}} \left[  	\sup_{f \in \mathcal{B}} |G_{n}[f]  - G^{\ast}_{n,q}[f]| \right] \leq & n^{-1/2} \sum_{j=0}^{n/q-1} E_{\mathbb{P}} \left[  \sup_{f \in \mathcal{B}} \sum_{l=1}^{q}  | f(X_{qj+l}) - f(X^{\ast}_{qj+l}) | \right] \\
		\leq &  n^{-1/2} \sum_{j=0}^{n/q-1} E_{\mathbb{P}} \left[ d_{q,\mathcal B} (U_{j+1},U^{\ast}_{j+1})   \right].
	\end{align*}
	Thus, by expression \ref{eqn:coupling.Hq.1},
	\begin{align*}
		E_{\mathbb{P}} \left[  	\sup_{f \in \mathcal{B}} |G_{n}[f]  - G^{\ast}_{n,q}[f]| \right] \leq  n^{-1/2} \sum_{j=0}^{n/q-1} \tau_{d_{q,\mathcal B}}(\sigma(U_1,\ldots,U_{j-1});U_{j+1}).
	\end{align*}
	
	It is clear that $\sigma(U_1,\ldots,U_{j-1}) \subseteq \sigma(X_1,\ldots,X_{jq-q}) = : \mathcal{M}_{1}^{jq-q}$, so $\tau_{d_{q,\mathcal B}}(\sigma(U_1,\ldots,U_{j-1});U_{j+1}) \le  \tau_{d_{q,\mathcal B}}(\mathcal{M}_{1}^{jq-q};(X_{jq+1},\ldots ,X_{jq+q} ))$ as $\tau_{d}$ is non-decreasing in the first argument. Thus, the desired result follows. 
\end{proof}

Following \cite{dedecker2006inequalities,dedecker2004coupling}, for any $q,i \in \mathbb{N}$, let
\begin{align}
 \tau_{q,\mathcal B}(i)   :  =\max_{1\leq l \leq q} \frac{1}{l} \sup \left\{  	\tau_{d_{l,\mathcal B}}( \mathcal{M}^{p}  ; X_{t_{1}},...,X_{t_{l}}) \colon p + i +1\leq t_{1} \leq \ldots \leq t_{l}   \right\}.
\end{align}
The pair $(\mathcal M^p,\; X_{j_1},\dots,X_{j_\ell})$ consists of the \emph{past} of the process up to time $p$,  $\mathcal M^p = \sigma(X_j : j \le p)$,  together with a \emph{future block} $(X_{j_1},\dots,X_{j_\ell})$ satisfying $j_1 \ge p+i+1$. The supremum is taken over all such blocks.

\begin{remark}
	In the text and throughout, when there is no risk of confusion, we use $\tau_{\mathcal B}(q)$ to denote $\tau_{q,\mathcal B}(q)$ for any $q \in \mathbb{N}$. 
\end{remark}

The following results is an immediate consequence of Lemma \ref{lem:tau-bdd.0} and the arguments in \cite{dedecker2006inequalities,dedecker2004coupling}. 

\begin{lemma}
	\label{lem:tau-bdd}
	Let $\mathcal{B} \subseteq cone \mathcal{F}$. For any $n \in \mathbb{N}$ and $q \in \mathcal{Q}_{n}$ such that $n/q \in \mathbb{N}$, 
	\begin{align*}
		\left \Vert \sup_{f \in \mathcal{B} } | G_{n}[f]  - G^{\ast}_{n,q}[f]     |  \right \Vert_{L^{1}(\mathbb{P})} \leq   \sqrt{n}  \tau_{q,\mathcal B}(q) .
	\end{align*}
	
\end{lemma}

\begin{proof}
	By Lemma \ref{lem:tau-bdd.0} and the fact that $\mathcal M_{1}^{jq-q} \subseteq \mathcal M^{jq-q}$, 
	\begin{align*}
		\left \Vert \sup_{f \in \mathcal{B} } | G_{n}[f]  - G^{\ast}_{n,q}[f]     |  \right \Vert_{L^{1}(\mathbb{P})} \leq  n^{-1/2} \sum_{j=0}^{n/q-1}   \tau_{d_{q,\mathcal B}}(\mathcal{M}^{jq-q};(X_{jq+1},\ldots ,X_{jq+q} )).
	\end{align*}
	
	Observe that
	\begin{align*}
		\tau_{d_{q,\mathcal B}}(\mathcal{M}^{jq-q};(X_{jq+1},\ldots ,X_{jq+q} ))  \leq & \max_{1\leq l \leq q} \sup \left\{  	\tau_{d_{l,\mathcal B}}(\mathcal{M}^{jq-q}; (X_{t_{1}},...,X_{t_{l}}))  \colon jq +1 \leq t_{1} \leq \ldots \leq t_{l}   \right\} \\
		\leq & q  \max_{1\leq l \leq q} \frac{1}{l} \sup \left\{  	\tau_{d_{l,\mathcal B}}(\mathcal{M}^{jq-q}; (X_{t_{1}},...,X_{t_{l}}))  \colon jq +1 \leq t_{1} \leq \ldots \leq t_{l}   \right\} \\
		= : & q \tau_{q,\mathcal B}(q).
	\end{align*}
	
	Therefore, 
	\begin{align*}
		\left \Vert \sup_{f \in \mathcal{B} } | G_{n}[f]  - G^{\ast}_{n,q}[f]     |  \right \Vert_{L^{1}(\mathbb{P})} \leq  n^{-1/2} \sum_{j=0}^{n/q-1} q \tau_{q,\mathcal B}(q) = \frac{q}{\sqrt{n}} \frac{n}{q}  \tau_{q,\mathcal B}(q)  = \sqrt{n}  \tau_{q,\mathcal B}(q) .
	\end{align*}
\end{proof}

\subsection{Useful Bounds for the $\tau$ coefficient} 


The following lemma contains some useful properties of the measure of dependence.

\begin{lemma}\label{lem:tau.properties}
	The following are true for any $q \in \mathbb{N}$,
	\begin{enumerate}
		\item For any $\mathcal{B} \subseteq \mathcal{B}' \subseteq  \mathcal{F}$, $ \tau_{\mathcal{B}}(q) \leq  \tau_{\mathcal{B}'}(q)$.
		\item For any $\mathcal{B} \subseteq  \mathcal{F}$, $\tau_{\mathcal{B}}(q) \leq \sup_{f \in  \mathcal{B}} ||f||_{L^{\infty}} \tau_{\bar{\mathcal{B}}}(q)$ where $\bar{\mathcal{B}} : = \{ f \in cone \mathcal{B} \colon ||f||_{L^{\infty}} \leq 1  \}$.
	\end{enumerate}
\end{lemma}

\begin{proof}[Proof of Lemma \ref{lem:tau.properties}]
	(1) It suffices to show that $\Lambda(\mathcal{B}) \subseteq \Lambda(\mathcal{B}')$. Take any $g \in  \Lambda(\mathcal{B})$, by definition of $\Lambda(\mathcal{B})$, for any $x,y$, $|g(x)-g(y)| \leq d_{\mathcal{B}}(x,y)$. Since $ d_{\mathcal{B}}(x,y) : = \sup_{h \in \mathcal{B}} |h(x) -h(y)| \leq \sup_{h \in \mathcal{B}'} |h(x) -h(y)| = : d_{\mathcal{B}'}(x,y)$ it follows that  $|g(x)-g(y)| \leq d_{\mathcal{B}'}(x,y)$, and thus $g \in  \Lambda(\mathcal{B}')$.
	
	\bigskip

	(2) If $\sup_{f \in  \mathcal{B}} ||f||_{L^{\infty}}  = \infty$ the inequality is trivial, so we proceed under the assumption that $B : = \sup_{f \in  \mathcal{B}} ||f||_{L^{\infty}} < \infty$.  Take any $g \in \Lambda(\mathcal{B})$ and observe that for any $x,y$, $|g(x)-g(y)| \leq \sup_{h \in \mathcal{B}} |h(x) -h(y)|$. Also observe that $ \sup_{h \in \mathcal{B}} |h(x) -h(y)| =  B \sup_{h \in \mathcal{B}}  |h(x)/B -h(y)/B| \leq B  \sup_{h' \in \bar{\mathcal{B}}}  |h'(x) -h'(y)|$, where the last inequality follows from the fact that $h/B \in cone \mathcal{B}$ and $||h/B||_{L^{\infty}} = B^{-1} ||h||_{L^{\infty}}  \leq 1$. Therefore, $g/B \in \Lambda(\bar{\mathcal{B}})$, and so $\tau_{\mathcal{B}}(q) = B \tau_{\mathcal{ \bar{B}}}(q)$
	as desired. 
\end{proof}


\subsubsection{Realtionship with Absolutely regular coefficient in \cite{DMR1995} (DMR).} 
\label{app:tau.bound}

For two $\sigma$-fields $\mathcal A,\mathcal B$, define
\[
\beta(\mathcal A,\mathcal B)
:=\frac12\,\mathbb E \left[  \Big\| \mathbb P(\,\cdot\,|\mathcal A)-\mathbb P(\,\cdot\,)\Big\|_{TV,\mathcal B} \right].
\]
In DMR one takes $\beta(q):=\beta\!\left(\mathcal M_{-\infty}^{0},\mathcal M_q^{+\infty}\right)$ for any $q \in \mathbb{N}$.

\begin{lemma}\label{lem:tau.bound.beta}
	Let $\mathcal{B} \subseteq cone \mathcal{F}$.  For any $q \in \mathbb{N}$, 
	\begin{align}
		\tau_{q,\mathcal B}(q) \leq  2 M 	\beta(q)
	\end{align}with $M:=\sup_{f\in\mathcal B}\|f\|_{L^{\infty}}$.
\end{lemma}

\begin{proof}[Proof of Lemma \ref{lem:tau.bound.beta}]
	For any $l \leq q$, let $X_{1:l}=(X_1,\dots,X_l)$ and define the block discrete metric	$d_0^{(l)}(x,y) =\sum_{\ell=1}^l \mathbf 1\{x_\ell\neq y_\ell\}$.  Observe that if $x_\ell=y_\ell$, then $\sup_{f\in\mathcal B}|f(x_\ell)-f(y_\ell)|=0$ and if $x_\ell\neq y_\ell$, then 
	\[
	\sup_{f \in \mathcal{B}} |f(x_\ell)-f(y_\ell)|
	\le \sup_{f \in \mathcal{B}} |f(x_\ell)|+|f(y_\ell)|
	\le 2\sup_{f \in \mathcal{B}} \|f\|_{L^{\infty}}
	\le 2M,
	\]
	with $	M=\sup_{f\in\mathcal B}\|f\|_{L^{\infty}}$. Thus, 	for any $x_{1:l} : = (x_{1} , \ldots , x_{l})$ and $y_{1:l} : = (y_{1} , \ldots , y_{l})$,
	\begin{align*}
		d_{l,\mathcal B}(x_{1:l},y_{1:l}) = \sum_{\ell=1}^{l} \sup_{f \in \mathcal{B}} |f(x_\ell)-f(y_\ell)| \leq 2 M \sum_{\ell=1}^{l} \mathbf 1\{x_\ell\neq y_\ell\} = 2M d_{0}^{(l)}(x_{1:l}, y_{1:l}) . 
	\end{align*}
	
	Thus, for any $g \in \Lambda_1(\mathbb{X}^{l},d_{l,\mathcal B})$, $|g( x_{1:l}  ) - g(y_{1:l}  ) | \leq 	d_{l,\mathcal B}(x_{1:l}  , y_{1:l}  )  \leq 2M d_{0}^{(l)}(x_{1:l}  , y_{1:l}  )$. So, $ \Lambda_1(\mathbb{X}^{l},d_{l,\mathcal B}) \subseteq  \Lambda_1(\mathbb{X}^{l},2M d_{0}^{(l)}) \subseteq 2 M \Lambda_1(\mathbb{X}^{l},d_{0}^{(l)})$. Hence,
	\begin{align*}
		\tau_{d_{l,\mathcal B}}(\mathcal{M}_{jq-q};(X_{t_{1}} , \ldots , X_{t_{l}}  )) \leq 2 M 	\tau_{d_{0}^{(l)}}(\mathcal{M}_{jq-q};(X_{t_{1}} , \ldots , X_{t_{l}}  )).
	\end{align*}
	Since any $g \in  \Lambda_1(\mathbb{X}^{l},d_{0}^{(l)})$ is such that $||g||_{L^{\infty}} \leq l/2$, it follows that
	\begin{align*}
		 \tau_{d_{l,\mathcal B}}(\mathcal{M}_{jq-q};(X_{t_{1}} , \ldots , X_{t_{l}}  )) \leq  & M l \left \Vert  \sup_{\|g\|_{L^{\infty}} \le 1} 	\left| 	\int g(x_{t_{1}} , \ldots , x_{t_{l}}  ) (P(d(x_{t_{1}} , \ldots , x_{t_{l}}) \mid \mathcal M_{jq-q} )- P (x_{t_{1}} , \ldots , x_{t_{l}}) )	\right | \right \Vert_{L^{1}(\mathbb{P})}\\
		 = & 2 M l \beta (\mathcal M_{jq-q},\sigma( X_{t_{1}} , \ldots , X_{t_{l}} ))
	\end{align*}
	for any $1\leq l \leq q$ and  $jq \leq t_{1} \leq ... \leq t_{l}$ in $\mathbb{N}$;  where the second line follows from the definition of $\beta$-coefficient and the results in \cite{dedecker2005new} and \cite{dedecker2006inequalities}.
	
	Because $\mathcal{M}_{-\infty}^{jq-q} \supseteq \mathcal M_{jq-q}$ and $\mathcal{M}_{jq}^{\infty} \supseteq \mathcal \sigma( X_{t_{1}} , \ldots , X_{t_{l}} )$ for any $jq \leq t_{1} \leq ... \leq t_{l}$ in $\mathbb{N}$; $\beta$ is non-decreasing in both arguments; and the process is stationary:
	\begin{align*}
		\frac{1}{l} \tau_{d_{l,\mathcal B}}(\mathcal{M}_{jq-q};(X_{t_{1}} , \ldots , X_{t_{l}}  )) \leq  2M 	\beta(q).
	\end{align*}
	Since this holds for any $1\leq l \leq q$ and  $jq \leq t_{1} \leq ... \leq t_{l}$ in $\mathbb{N}$, one can take the maximum in the LHS and obtain
	\begin{align*}
		\tau_{q,\mathcal B}(q) \leq  2M 	\beta(q).
	\end{align*}
\end{proof}

\subsubsection{Lipschitz Classes}
\label{app:Lipschitz}

 	Suppose $\mathcal B$ is a set of measurable
functions $f:\mathsf X\to\mathbb R$ such that
\begin{equation}\label{eq:B-rho-domination}
	|f(x)-f(y)| \le \rho(x,y),
	\qquad \forall f\in\mathcal B,\ \forall x,y\in\mathsf X,
\end{equation}
for a metric $\rho$ over $\mathsf X$. 

For any $l \in \mathbb{N}$, define the block pseudo-metric on $\mathsf X^l$ by $d_{l,\mathcal B}(x_{1:l},y_{1:l})
:=
\sum_{m=1}^{l}\ \sup_{f\in\mathcal B}\,|f(x_m)-f(y_m)|$, 
and define the $\rho$-sum metric on $\mathsf X^l$ by
$\rho^{(l)}(x_{1:l},y_{1:l})
:=
\sum_{m=1}^{l}\rho(x_m,y_m)$. 

It follows that $ \Lambda_1(\mathsf X^{l},d_{l,\mathcal B}) \subseteq  \Lambda_1(\mathsf X^{l},\rho^{(l)}) $. Thus,  for any $p \in \mathbb{N}$ and any block $(X_{t_{1}},...,X_{t_{l}})$,
\begin{align*}
	\tau_{d_{l,\mathcal B}}( \mathcal{M}^{-q}_{-\infty}  ; X_{t_{1}},...,X_{t_{l}}) \leq 	\tau_{\rho^{(l)}}( \mathcal{M}^{-q}_{-\infty}  ; X_{t_{1}},...,X_{t_{l}}).
\end{align*}
This readily implies
\begin{align*}
	\tau_{q,\mathcal B}(q) \leq 	\tau_{\rho,q}(q) 	
\end{align*}
where, abusing notation, 
\begin{align*}
	\tau_{\rho,q}(q)   :  =\max_{1\leq l \leq q} \frac{1}{l} \sup \left\{  	\tau_{\rho^{(l)}}( \mathcal{M}^{-q}_{-\infty}  ; X_{t_{1}},...,X_{t_{l}}) \colon 1\leq t_{1} \leq \ldots \leq t_{l}   \right\}.
\end{align*}
And also, 
by stationarity, 
\begin{align}
	\tau_{d_{q,\mathcal B}}(\mathcal{M}^{jq-q}_{1};(X_{jq+1},\ldots ,X_{jq+q} )) \leq 	\tau_{d_{q,\mathcal B}}(\mathcal{M}^{jq-q}_{-\infty};(X_{jq+1},\ldots ,X_{jq+q} ))   \leq   \tau_{\rho^{(q)} }  (\mathcal{M}^{-q}_{-\infty};(X_{1},\ldots ,X_{q} )).
\end{align}

This result applied to $\mathsf X = \mathbb X$ and Lemma \ref{lem:tau-bdd.0} readily imply
\begin{lemma}
	\label{lem:tau-bdd.Lip}
	Let $\mathcal{B} \subseteq cone \mathcal{F}$ such that expression \ref{eq:B-rho-domination} holds. For any $n \in \mathbb{N}$ and $q \in \mathcal{Q}_{n}$ such that $n/q \in \mathbb{N}$, 
	\begin{align*}
		\left \Vert \sup_{f \in \mathcal{B} } | G_{n}[f]  - G^{\ast}_{n,q}[f]     |  \right \Vert_{L^{1}(\mathbb{P})} \leq  \sqrt{n}  \frac{1}{q}  \tau_{\rho^{(q)} }  (\mathcal{M}^{-q}_{-\infty};(X_{1},\ldots ,X_{q} )).
	\end{align*}
And 
	\begin{align*}
	\left \Vert \sup_{f \in \mathcal{B} } | G_{n}[f]  - G^{\ast}_{n,q}[f]     |  \right \Vert_{L^{1}(\mathbb{P})} \leq  \sqrt{n}   \tau_{\rho,q}  (q).
\end{align*}
\end{lemma}

\section{Supplementary Lemmas}
\label{app:supp.lemmas} \label{app:Qn}

%
%

Recall that $\alpha^{-1}$ denotes the generalized inverse of $\alpha$ defined in expression \ref{eqn:alpha.coeff}, $q \mapsto \bar{\tau}(q) : = \tau(q)/q$, and $(\mathbf{q}_{n}(k))_{k}$ is defined in expression \ref{eqn:qk} with $\mathcal{Q}_{n}$ being the set of integers less or equal than $n$.

The next result presents an useful and known (cf. \cite{DMR1995}) bound for $\mu_{q}$ --- defined in expression \ref{eqn:muq}. 

\begin{lemma}\label{lem:mu.bound}
	For any $q \in \mathbb{N}$,
	\begin{align*}
		\mu_{q}(u) \leq \min\{ \alpha^{-1}(u) , q \} + 1,~\forall u \in [0,1].
	\end{align*}
\end{lemma}

\begin{proof}
	 Observe that $u \mapsto \mu_{q}(u) = \sum_{l=0}^{q} 1\{ u \leq \alpha(l) \}$, and so, for any $u \in [0,1]$,  $\mu_{q}(u) = \min\{ j(u), 1+q\}$ where  $j(u) = \min \{ s \in \mathbb{N}_{0} \mid u > \alpha(s)  \}$. Hence, it suffices to show that $\alpha^{-1}(u) \geq j(u)-1$. By definition of $j(.)$, $u \leq \alpha(j(u)-1)$ and by definition of $\alpha^{-1}(u)$, $\alpha(\alpha^{-1}(u)) \leq u$; thus $\alpha(\alpha^{-1}(u)) \leq \alpha(j(u)-1)$. Since $\alpha$ is non-increasing, this implies that $\alpha^{-1}(u) \geq j(u)-1$. Thus, the desired results holds. 	
\end{proof}

The next result presents an useful bound for $\mathbf{q}_{n}(0)$. 

\begin{lemma}\label{lem:q0.bound}
	For any $n \in \mathbb{N}$, $\mathbf{q}_{n}(0) \leq \bar{\tau}^{-1}(1/n) + 1$.
\end{lemma}

\begin{proof}
%
%

	By definition $\mathbf{q}_{n}(0) = \min \{ q \in \mathcal{Q}_{n} \colon \bar{\tau}(q) \leq 1/n \}$. Thus, either $\mathbf{q}_{n}(0) = 1$, or $ \bar{\tau}(\mathbf{q}_{n}(0)-1) > 1/n$. This last inequality and the fact that $\bar{\tau}$ is non-increasing readily imply that $\mathbf{q}_{n}(0) \leq \bar{\tau}^{-1}(1/n) + 1$.
\end{proof}

The next lemma relates summability of $\tau$ with a convergence rate for $(\mathbf{q}_{n}(0))_{n}$.

\begin{lemma}\label{lem:tau.bounded.qn.rate}
	If $\sum_{m=1}^\infty \tau(m)<\infty$ and $	x(n):=\min\{x \le n:\ \tau(x)\le x/n\}$, then $\limsup_{n\to\infty} \frac{x(n)}{\sqrt{n}} = 0$.
\end{lemma}

\begin{proof}
	Suppose, by contradiction, that \(\limsup_{n\to\infty} x(n)/\sqrt{n}\ge \varepsilon>0\).
	Then there exists an increasing  and diverging sequence $(n_k)_{k}$ such that
	\[
	x(n_k)\ \ge\ \varepsilon\sqrt{n_k}\qquad (k\ge1).
	\]
	Set \(M_k:=\lfloor \varepsilon\sqrt{n_k}\rfloor\), so \(M_k\uparrow\infty\) and \(M_k\le x(n_k)\).
	By the minimality of \(x(n_k)\), for every \(x\le M_k\) we must have $\tau(x)\ >\ \frac{x}{n_k}$. 
	Hence, for \(k\ge1\),
	\[
	\sum_{x=M_{k-1}+1}^{M_k}\tau(x)
	\ \ge\
	\frac{1}{n_k}\sum_{x=M_{k-1}+1}^{M_k} x
	\ =\
	\frac{1}{2n_k}\Big(M_k(M_k+1)-M_{k-1}(M_{k-1}+1)\Big).
	\]
	Since \(M_k=\varepsilon\sqrt{n_k}+O(1)\), there exists \(k_0\) and a constant \(c_0\in(0,\varepsilon^2]\)
	such that for all \(k\ge k_0\),
	\[
	\frac{1}{2n_k}\Big(M_k^2-M_{k-1}^2\Big)\ \ge\
	\frac{c_0}{2}\Big(1-\frac{n_{k-1}}{n_k}\Big).
	\]
	Therefore, summing over disjoint blocks,
	\[
	\sum_{m=1}^\infty \tau(m)
	\ \ge\
	\sum_{k=k_0}^\infty \sum_{x=M_{k-1}+1}^{M_k}\tau(x)
	\ \ge\
	\frac{c_0}{2}\sum_{k=k_0}^\infty\Big(1-\frac{n_{k-1}}{n_k}\Big).
	\]
	The series on the right diverges because for \(0<r\le1\) one has \(-\log r\ge 1-r\), hence
	\[
	\sum_{k=k_0}^K \Big(1-\frac{n_{k-1}}{n_k}\Big)
	\ \ge\
	\sum_{k=k_0}^K -\log\!\Big(\frac{n_{k-1}}{n_k}\Big)
	\ =\
	-\log\!\Big(\frac{n_{k_0-1}}{n_K}\Big)
	\ \xrightarrow[K\to\infty]{}\ \infty,
	\]
	since \(n_K\to\infty\). This contradicts \(\sum_{m=1}^\infty \tau(m)<\infty\).
\end{proof}

\section{Proofs of  Lemmas in Section \ref{sec:proof.main}}
\label{app:main.proof}

\begin{proof}[Proof of Lemma \ref{lem:bere}]
	Throughout, we use $J = J(q,n) = n/q$, which is an integer since $n/q$ is assumed to be an integer in this lemma.	Let $f \mapsto \delta_{j}(\omega^{\ast})f := q^{-1/2} \sum_{i=qj+1}^{qj+q} \{f(X^{\ast}_{i}) - E_{P}[f(X^{\ast})]\}$ for any $j = 0,1,2,...$. Then, it follows that for any $g \in \mathcal{F}$, 
	\begin{align*}
		G^{\ast}_{n,q}[g] =& n^{-1/2} \sqrt{q} \sum_{j=0}^{J-1} \delta_{j}(\omega^{\ast})g
	\end{align*}
	where $\delta_{j}(\omega^{\ast})g$ is measurable with respect to $U^{\ast}_{j}(q)=(X^{\ast}_{qj+1},...,X^{\ast}_{qj+q})$. 
	
	Note that, for any sequence $(a_{j})_{j} $, $\sum_{j=0}^{J-1} a_{j} = \sum_{j=0}^{[(J-1)/2]}  a_{2j} + \sum_{j=0}^{[J/2]-1} a_{2j+1}$, then \footnote{Here, $[a] = \max\{ n \in \mathbb{N} : n \leq a  \}$. } 					
	\begin{align}\label{eqn:bere-1}
		G^{\ast}_{n,q}[g] = n^{-1/2} \sqrt{q} \left(   \sum_{j=0}^{[(J-1)/2]}  \delta_{2j}(\omega^{\ast}) g  + \sum_{j=0}^{[J/2]-1}  \delta_{2j+1}(\omega^{\ast}) g   \right).
	\end{align}
	Hence, for any $t>0$, \begin{align}
		P \left(	G^{\ast}_{n,q}[g] \geq t \right) \leq P \left(	n^{-1/2} \sqrt{q}   \sum_{j=0}^{[(J-1)/2]}  \delta_{2j}(\omega^{\ast}) g    \geq 0.5 t \right) 
		+P  \left(	n^{-1/2} \sqrt{q}  \sum_{j=0}^{[J/2]-1}  \delta_{2j+1}(\omega^{\ast})g  \geq 0.5 t \right) \label{eqn:bere-2}.
	\end{align}
	
	Thus, it suffices to bound each term in the right hand side (RHS) separately. By the results in Appendix \ref{app:coupling}, $(U^{\ast}_{2j}(q))_{j \geq 0}$ are independent (so are $(U^{\ast}_{2j+1}(q))_{j \geq 0}$). Moreover, they  have the same distribution as $(U_{2j}(q))_{j\geq 0 }  $ (also $(U^{\ast}_{2j+1}(q))_{j \geq 0}$ has the same distribution as $(U_{2j+1}(q))_{j \geq 0}$). This means that both sums are comprised on IID terms. 
	
	It is easy to see that for any $g \in \mathcal{F}$, $||\delta_{j}(\omega^{\ast})g||_{L^{\infty}} \leq \sqrt{q} 2 ||g||_{L^{\infty}}$. Moreover,  for any $g \in \mathcal{F}$,	
	\begin{align} \notag 
		E_{P}[(\delta_{0}(\omega^{\ast}) g )^{2}] 	= & Var(g(X_{0})) + 2 \sum_{k=1}^{q} (1-k/q) Cov(g(X_{0}),g(X_{k}))\\ \notag
		\leq & 4 \sum_{k=0}^{q} (1-k/q) \int_{0}^{\alpha(g(X_{k}),g(X_{0}))} Q^{2}_{g}(u) du \\ \notag
	\leq & 4 \sum_{k=0}^{q} (1-k/q) \int_{0}^{\alpha(k)} Q^{2}_{g}(u) du \\ \notag
		\leq & 4 \int_{0}^{1} \left(\sum_{k=0}^{q} 1 \{ u \leq  \alpha(k)   \}    \right) Q^{2}_{g}(u) du\\ \label{eqn:variance.bound}
		\leq & 4 \int_{0}^{1} \mu_{q}(u) Q^{2}_{g}(u) du = 4 ||g||_{2,q}^{2}
	\end{align}
where the second line follows by \cite{rio2017asymptotic} Theorem 1.1 where $\alpha(g(X_{k}),g(X_{0}))$ is defined in  \cite{rio2017asymptotic} in 1.18a and by \cite{dedecker2005new} equals $\sup_{f : \mathbb{R} \rightarrow \mathbb{R} \colon f \in BV_{1}} \left|  \int f(g(X)) P(dx \mid \mathcal{M}_{-k}) -  \int f(g(X)) P_{X}(dx)  \right|$. The third line follows because this quantity is bounded by $\alpha(k)$ defined in expression \ref{eqn:alpha.coeff}.


	

 Hence,	by Bernstein inequality (see \cite{VdV-W1996} Lemma 2.2.9), it follows that
	\begin{align*}
		P  \left(  \sqrt{\frac{q}{n}} \sum_{j=0}^{[(J-1)/2]}  \delta_{2j}(\omega^{\ast})g   \geq t         \right) \leq &
		\exp \left\{ - \frac{0.5 q^{-1}n t ^{2}}{\sum_{j=0}^{[(J-1)/2]} E_{P}[(\delta_{2j}(\omega^{\ast}) g )^{2}]  + \frac{2}{3} \sqrt{n} t ||g||_{L^{\infty}}    }  \right\} \\
		\leq & \exp \left\{ - \frac{0.5 n q^{-1} t ^{2}}{(0.5J+1) 4||g||_{2,q}^{2}  + \frac{2}{3} \sqrt{n} t   ||g||_{L^{\infty}}    }  \right\}\\
		\leq & \exp \left\{ - \frac{0.5  t ^{2}}{4||g||_{2,q}^{2}  + \frac{2}{3} t \frac{q}{\sqrt{n}} ||g||_{L^{\infty}}    }  \right\}
	\end{align*}
	where the second line follows from the fact that $E_{P}[(\delta_{2j}(\omega^{\ast})f )^{2}] = E_{P}[(\delta_{0}(\omega^{\ast}) f )^{2}]$ (stationarity) and the bound in \ref{eqn:variance.bound}, the third line follows from simple algebra and the fact that $J = \frac{n}{q}$ and $n/q \geq 1$.
	
	By applying the same calculations to $P \left(  \sqrt{\frac{q}{n}} \sum_{j=0}^{[J/2]-1}  \delta_{2j+1}(\omega^{\ast})g   \geq t         \right) $ it follows that 
	\begin{align}
	P \left(	G^{\ast}_{n,q}[g] \geq t \right) \leq 2 \exp \left\{ - \frac{t ^{2}}{ 32  ||g||_{2,q}^{2}  + \frac{16}{3} t \frac{q}{\sqrt{n}} ||g||_{L^{\infty}}    }  \right\}.
	\end{align}

	By letting $t  = \sqrt{32 ||g||^{2}_{2,q} u^{2} } + \frac{16}{3} \frac{q}{\sqrt{n}} ||g||_{L^{\infty}} u^{2} $ for any $u > 0$, our result follows as
	\begin{align}
	P  \left(	G^{\ast}_{n,q}[g] \geq \sqrt{32 ||g||^{2}_{2,q} u^{2} } + \frac{16}{3} \frac{q}{\sqrt{n}} ||g||_{L^{\infty}} u^{2} \right) \leq 2 e^{-u^{2}},
	\end{align}
To show this claim, observe that 
	\begin{align*}
		\frac{t ^{2}}{32||g||_{2,q}^{2}  + \frac{16}{3} t \frac{q}{\sqrt{n}} ||g||_{L^{\infty}}    } = &  u^{2} \left( \frac{32 ||g||^{2}_{2,q} + \left( \frac{16}{3} \frac{q}{\sqrt{n}} ||g||_{L^{\infty}} \right)^{2} u^{2} + 2 \frac{16}{3} \frac{q}{\sqrt{n}} ||g||_{L^{\infty}} \sqrt{32 ||g||^{2}_{2,q}} u  }{32 ||g||^{2}_{2,q} +  \left(\frac{16}{3} \frac{q}{\sqrt{n}} ||g||_{L^{\infty}}  \right)^{2} u^{2}    + \frac{16}{3} \frac{q}{\sqrt{n}} ||g||_{L^{\infty}} \sqrt{32 ||g||^{2}_{2,q}} u     }  \right)\\
		\geq & u^{2}.
	\end{align*}
\end{proof}

\begin{proof}[Proof of Lemma \ref{lem:bound.pr.talagrand}]
	We show that $$	P \left(  \sup_{f \in \mathcal{F}} \sum_{l=1}^{\infty}  \sqrt{\frac{ \underline{n}(l) }{n}}| G^{\ast}_{\underline{n}(l),\mathbf{q}(l)}[\Delta_{l} f] | \geq v   \sup_{f \in \mathcal{F}}  \sum_{l=1}^{\infty}   \sqrt{\frac{ \underline{n}(l) }{n}} \left( 6 ||\Delta_{l} f||_{2,\mathbf{q}(l)} \ell(l)   + \frac{16}{3} \frac{\mathbf{q}(l)}{\sqrt{\underline{n}(l)}} ||\Delta_{l} f||_{L^{\infty}} (\ell(l)  )^{2}    \right)   \right) \leq \mathbb{L}_{0} e^{-0.5 v }.$$ To do this, we observe that the LHS is bounded above by $P( \cup_{l \in \mathbb{N}} \cup_{f \in \mathcal{F}} 	\Omega(\sqrt{v},l,\Delta_{l} f )  )$, where for any $t \geq 4$, any $l \in \mathbb{N}$, and any $g \in \Delta_{l}\mathcal{F}$,
	\begin{align*}
		\Omega(t,l,g) :  = \left\{ |G^{\ast}_{\underline{n}(l),\mathbf{q}(l)}[g]  | \geq t^{2} \left( 6 ||g||_{2,\mathbf{q}(l)} \ell(l)   + \frac{16}{3} \frac{\mathbf{q}(l)}{\sqrt{\underline{n}(l)}} ||g||_{L^{\infty}} (\ell(l)  )^{2}  \right)  \right\}.
	\end{align*}
	Since $ t\geq 1$, the RHS of the inequality is bounded below by $ t \ell(l)   6 ||g||_{2,\mathbf{q}(l)}   + (t \ell(l) )^{2} \frac{16}{3} \frac{\mathbf{q}(l)}{\underline{n}(l)} ||g||_{L^{\infty}}   $. Thus, by Lemma \ref{lem:bere} --- applied to $n = \underline{n}(l)$,  $g = \Delta_{l}f$, $q =\mathbf{q}(l)$, and $u = \ell(l) t$ --- it follows that $P(	\Omega(t,l,g) ) \leq 2 e^{-(\ell(l) t)^{2}}$. Therefore applying the union bound for any $t \geq 2$,
	\begin{align*}
		P( \cup_{l \in \mathbb{N}} \cup_{g \in \Delta_{l} \mathcal{F}} 	\Omega(t,l,g)  )  \leq 	\sum_{l=1}^{\infty} \sum_{g \in \Delta_{l} \mathcal{F}} 2 e^{-2^{(\ell(l)t)^{2}}}.
	\end{align*}
	
	The cardinality of the set $\Delta_{l} \mathcal{F}$ is bounded by $ 2^{2d 2^{l/c}}$. To show this claim observe that $card \Delta_{l} \mathcal{F}  \leq  card \mathcal{T}_{l} \times card \mathcal{T}_{l-1} $, and $card \mathcal{T}_{l} \leq 2^{d2^{l/c}}$. Thus $ card \Delta_{l} \mathcal{F}  \leq 2^{d(2^{l/c} + 2^{(l-1)/c})} \leq 2^{d 2^{l/c} (1 + 2^{-1/c})} \leq 2^{2d 2^{l/c}}$.
	
	Since $(\ell(l) t)^{2} = 0.5 \ell(l)^{2}  t^{2} + 0.5  \ell(l)^{2} t^{2} \geq  2 \ell(l)^{2} + 0.5 t^{2}$ because $t^{2} \geq 4$ and $\ell(l) \geq 1$, the previous display implies 
	\begin{align*}
		P ( \cup_{l \in \mathbb{N}} \cup_{g \in  \Delta_{l} \mathcal{F}} 	\Omega(t,l, g )  ) \leq 2	 e^{-0.5 t^{2}} \sum_{l=0}^{\infty} 2^{  2d 2^{l/c}  }  e^{- 2\ell(l)^{2} }.
	\end{align*}
	
	By assumption  $d 2^{l/c} \leq \ell(l)^{2} $ for any $l \in \mathbb{N}$, hence $ \sum_{l=0}^{\infty} 2^{  2d 2^{l/c}  }  e^{- 2\ell(l)^{2} } < \infty$ and the desired result follows by setting $t^{2} = v$.
\end{proof}

\begin{proof}[Proof of Lemma \ref{lem:bound.E.talagrand}]
	
	By the proof of Lemma \ref{lem:bound.pr.talagrand},
	\begin{align}
		P \left(  \sup_{f \in \mathcal{F}} \left\{  \sum_{l=1}^{\infty} \sqrt{\frac{\underline{n}(l)}{n}}| G^{\ast}_{\underline{n}(l),\mathbf{q}(l)}[\Delta_{l} f] |   \right\} \leq v A   \right) \geq 1 - \left( 2 \sum_{l=0}^{\infty} 2^{  2d 2^{l/c}  }  e^{- 2\ell(l)^{2} } \right)  e^{-0.5 v}
	\end{align}
	where $A  : =   \sup_{f \in \mathcal{F}}  \sum_{l=1}^{\infty}   \sqrt{\frac{\underline{n}(l)}{n}}  \left( 6 ||\Delta_{l} f||_{2,\mathbf{q}(l)} \ell(l)   + ||\Delta_{l} f ||_{L^{\infty}} \left( \frac{16}{3} \frac{\mathbf{q}(l)}{\sqrt{\underline{n}(l)}}  ( \ell(l)  )^{2}    \right)   \right)  $.
	
	Observe that for any random variable $Z$ such that $P(|Z| \geq u A) \leq C e^{-0.5u}$ for any $u>0$ and constants $A,B,C$,
	\begin{align*}
		E_{P} \left[   |Z|   \right] =  \int_{0}^{\infty} P \left( |Z|  \geq u  \right) du  =  A \int_{0}^{\infty} P \left( | Z |  \geq A t   \right) dt   = A C \int_{0}^{\infty} e^{-0.5t } dt.
	\end{align*}

   Therefore, by taking $Z: = \sup_{f \in \mathcal{F}} \sum_{l=1}^{\infty}  \sqrt{\frac{\underline{n}(l)}{n}} G^{\ast}_{\underline{n}(l),\mathbf{q}(l)}[\Delta_{l} f] |    $, and $C := 2 \sum_{l=0}^{\infty} 2^{  2d 2^{l/c}  }  e^{- 2\ell(l)^{2} }$, the desired result follows.

\end{proof}

\section{Majorization of $\gamma_{\alpha,\beta}$ in terms of Metric Entropy}
\label{app:bound.Dudley}

For any metric space $(T,d)$, let $N(\epsilon,T,d)$ be the packing number associated with a radius $\epsilon>0$ defined as 
\begin{align*}
	N(\epsilon,T,d) = \inf\{ |F| \colon F \subseteq T~finite,~and~\sup_{t \in T} d(t,F) \leq \epsilon   \}.
\end{align*}
Also, for any $\beta>0$, let the entropy number for any $l \in \mathbb{N}_{0}$ be
\begin{align*}
	e_{l,\beta}(T,d) : = \inf_{S \colon S \subseteq T~with~card S \leq 2^{2^{l/\beta}} }   \sup_{t \in T} d(t,S).
\end{align*}

The following lemma is based on the calculations in \cite{talagrand2014} p. 21-22 and is here mostly for completeness. 
\begin{lemma}\label{lem:bound.Dudley}
	For any $(\alpha ,\beta) \in \mathbb{R}^{2}_{++}$ 
	\begin{align*}
	 	\gamma_{\alpha,\beta}(T,d)	\leq  \sum_{l=0}^{\infty} 2^{l/\alpha} e_{l,\beta}(T,d).
	\end{align*}
	And if $\beta \leq 0.5 \alpha$, then
	\begin{align*}
		\gamma_{\alpha,\beta}(T,d) \leq (1 - 2^{-1/(2\beta)})^{-1} \int_{0}^{Diam(T,d)} \sqrt{ \log N(\epsilon,T,d) }  d\epsilon.
	\end{align*}
\end{lemma}

\begin{proof}
	\textbf{Part 1.}  Follows directly from the derivations in p.20-21 in \cite{talagrand2014}. 
	
\bigskip
	
	\textbf{Part 2.} For any sequence decreasing sequence, $(a_{m})_{m}$ converging to $0$ and such that $a_{0} \leq Diam(T,d)$ it follows that 
	\begin{align}\label{eqn:Dudley.1}
		\int_{0}^{Diam(T,d)} \sqrt{ \log N(\epsilon,T,d) }  d\epsilon = \sum_{m=0}^{\infty}  \int_{a_{m+1}}^{a_{m}} \sqrt{ \log N(\epsilon,T,d) }  d\epsilon \geq \sum_{m=0}^{\infty}  \sqrt{ \log N(a_{m},T,d) }   (a_{m} - a_{m+1})
	\end{align}	
	where the inequality follows from the fact that $\epsilon \mapsto N(\epsilon,T,d)$ is non-increasing.
	
	Simplifying notation, for any $l \in \mathbb{N}_{0}$, let $e_{l} = 	e_{l,\beta}(T,d) $. Clearly, $e_{0} = \inf_{t_{0} \in T}  \sup_{t \in T} d(t,t_{0}) \leq Diam(T,d)$ and $2 e_{0} \geq Diam(T,d)$. Thus, in display \ref{eqn:Dudley.1}, we can $a_{m}=e_{m}$ for all $m \in \mathbb{N}_{0}$ and obtain
	\begin{align}\label{eqn:Dudley.2}
		\int_{0}^{Diam(T,d)} \sqrt{ \log N(\epsilon,T,d) }  d\epsilon  \geq \sum_{m=0}^{\infty}  \sqrt{ \log N(e_{m},T,d) }   (e_{m} - e_{m+1}).
	\end{align}

	We claim that $ N(e_{m},T,d)   \geq 2^{2^{m/\beta}}$. Suppose not. Then there exists a finite set $F$ such that $\sup_{t \in T} d(f,F) \leq e_{m}$ with $card F \leq 2^{2^{m/\beta}}-1$. But then, one can add one element to $F$ to form a set $S$ such that $card S \leq 2^{2^{m/\beta}}$ and  $\sup_{t \in T} d(f,S) < e_{m}$, but this violates the definition of $e_{m}$.
	
	Under this claim, $ \sqrt{ \log N(a_{m},T,d) }  \geq 2^{m/(2\beta)}$. Therefore by expression \ref{eqn:Dudley.2} it follows that
	\begin{align*}
		\int_{0}^{Diam(T,d)} \sqrt{ \log N(\epsilon,T,d) }  d\epsilon  & \geq  \sum_{m=0}^{\infty}  2^{m/(2\beta)}  (e_{m} - e_{m+1}) \\
		& = e_{0}  + e_{1}(  2^{1/(2\beta)}  - 2^{0/(2\beta)}  ) + e_{2}(  2^{2/(2\beta)}  - 2^{1/(2\beta)}  )  + .... \\
		& \geq (1 - 2^{-1/(2\beta)}) \sum_{l=0}^{\infty} 2^{l/(2\beta)} e_{l} \\
		& \geq  (1 - 2^{-1/(2\beta)}) \sum_{l=0}^{\infty} 2^{l/\alpha} e_{l} 
	\end{align*}	
	where the last line follows from the assumption that $2 \beta \leq \alpha$. The desired result follows from part 1.
\end{proof}

\begin{remark}[Bound of $\gamma_{1,b}$ under $L^{\infty}$ norm]\label{rem:bound.gamma1}
	There are several techniques to establish bounds for $\gamma_{1,b}(\mathcal{F},||.||_{L^{\infty}})$. The first expression in Lemma \ref{lem:bound.Dudley} provides one, but there are alternatives based on more sophisticated methods. For instance, Propositions 2.3 and Corollary 2.7 in \cite{VanHandel2018b} (see also \cite{VanHandel2018}) provide a contraction principle to establish bounds for $\gamma_{1,b}$ for $b=1$.  In particular, they provide sharper bounds than Lemma \ref{lem:bound.Dudley} when $T$ satisfies certain convexity restrictions.  
	
	Of particular interest for this paper and empirical process theory in general is the case where $\mathcal{F}$ is a function class with certain smoothness. Here we cover several examples that satisfy finiteness of $\gamma_{1,1}(\mathcal{F},||.||_{L^{\infty}})$ assuming there is enough smoothness. In what follows, let	$D \subseteq \mathbb{R}^{d}$ bounded with smooth boundary:
	\begin{itemize}
		\item $\mathcal{F} = B^{\alpha}_{p,q}(D)$ is the Besov space over $D$ with parameters $(p,q,\alpha)$ (see \cite{EdmundsTriebel1996} Ch. 2.2 for the formal definition). By \cite{EdmundsTriebel1996} Theorem 3.3.2 (with $s_{2} = 0, s_{1}=\alpha, p_{2} = \infty,p_{1}=p, q_{2}=\infty,q_{1}=q$) under the condition $\alpha/d > 1/p$,\footnote{Observe that in the \cite{EdmundsTriebel1996} the entropy numbers are defined using a cardinality of $2^{k}$ as opposed to $2^{2^{k}}$ as we do here.} 
		\begin{align*}
			e_{l}(\mathcal{F},||.||_{L^{\infty}}) \leq \mathbb{L} 2^{-l\alpha/d}.
		\end{align*}
		Assuming $\alpha/d > 1$ this readily implies that
		\begin{align*}
			\gamma_{1,1}(\mathcal{F},||.||_{L^{\infty}})	\leq  \mathbb{L} \sum_{l=0}^{\infty} 2^{l-l \alpha/d}  < \infty. 
		\end{align*} 
		 
		\item $\mathcal{F} = F^{\alpha}_{p,q}(D)$ is the Triebel-Lizorkin space over $D$ with parameters $(p,q,\alpha)$ (see \cite{EdmundsTriebel1996} Ch. 2.2 for the formal definition). The same result as for Besov spaces holds here (see \cite{EdmundsTriebel1996} Theorem 3.3.2). 
		
		\item $\mathcal{F} = \mathbb{W}^{\alpha}_{p}(D)$ is the Sobolev space over $D$ with parameters $(p,\alpha)$. Since $\mathbb{W}^{\alpha}_{p}(D)=B^{\alpha}_{p,p}(D)$, we also obtain $	\gamma_{1,1}(\mathcal{F},||.||_{L^{\infty}}) < \infty$ in this case under the assumption $\alpha/d  > 1$. 
		
		\item $\mathcal{F} = \mathcal{C}^{\alpha}(D)$ is the H\"{o}lder-Zygmund space. Since $\mathcal{C}^{\alpha}(D) = B^{\alpha}_{\infty,\infty}(D)$, we also obtain $	\gamma_{1,1}(\mathcal{F},||.||_{L^{\infty}}) < \infty$ in this case under the assumption $\alpha/d  > 1$.\footnote{These results also hold for $D=[0,1]^{d}$ or any bounded domain with non-empty interior.} 
	\end{itemize}
	We refer the reader to \cite{EdmundsTriebel1996} for bounds on entropy numbers for other function classes; in particular for ``weighted" versions of the aforementioned spaces.$\triangle$
\end{remark}

%
%
%
%
%
%
%
%
%
%
%
%
%
%

	\end{document}